\documentclass[11pt]{article}
\usepackage{amsmath,amsthm,amsfonts,verbatim,amssymb}

 \DeclareMathOperator{\lca}{lca}
 \DeclareMathOperator{\diam}{diam}
 \DeclareMathOperator{\Lip}{Lip}
 \DeclareMathOperator{\coLip}{coLip}
 \DeclareMathOperator{\geo}{\text{geo}}

\newcommand{\E}{\mathbb{E}}
\newcommand{\U}{\mathcal{U}}
\newcommand{\HH}{\mathcal{H}}
\newcommand{\V}{\mathcal{V}}
\newcommand{\W}{\mathcal{W}}
\newcommand{\M}{\mathcal{M}}
\newcommand{\N}{\mathbb{N}}
\newcommand{\A}{\mathcal{A}}
\newcommand{\C}{\mathcal{C}}
\newcommand{\R}{\mathbb{R}}
\newcommand{\Q}{\mathcal{Q}}
\newcommand{\QS}{\mathcal{QS}}
\newcommand{\SQ}{\mathcal{SQ}}
\newcommand{\e}{\varepsilon}

\newcommand{\comp}{\text{comp}}

\theoremstyle{plain}
 \newtheorem{lemma}{Lemma}[section]
 \newtheorem{theorem}[lemma]{Theorem}
 \newtheorem{corollary}[lemma]{Corollary}
 \newtheorem{claim}[lemma]{Claim}
 \newtheorem{proposition}[lemma]{Proposition}

 \theoremstyle{definition}
 \newtheorem{definition}[lemma]{Definition}

 \newtheorem{remark}[lemma]{Remark}

\begin{document}

\title{Euclidean Quotients of Finite Metric Spaces}

\author{
    Manor Mendel%
\thanks{Supported in part by a grant from the Israeli Science Foundation
(195/02) and by the Landau Center.}\\
    School of Computer Science, The Hebrew University, Jerusalem. 
    \and
    Assaf Naor \\
    Microsoft Research, Redmond, Washington \\
     }

\date{}

\maketitle

\begin{abstract}
This paper is devoted to the study of quotients of finite metric
spaces. The basic type of question we ask is: Given a finite
metric space $M$ and $\alpha\ge 1$, what is the largest quotient
of (a subset of) $M$ which well embeds into Hilbert space. We
obtain asymptotically tight bounds for these questions, and prove
that they exhibit phase transitions. We also study the analogous
problem for embedings into $\ell_p$, and the particular case of
the hypercube.
\end{abstract}

\bigskip { \em `` Our approach to general metric spaces
bears the undeniable imprint of early exposure to Euclidean
geometry. We just love spaces sharing a common feature with
$\R^n$." \flushright Misha Gromov.

}

\section{Introduction}\label{introduction}

A classical theorem due to A. Dvoretzky states that for every
$n$-dimensional normed space $X$ and every $\e>0$ there is a
linear subspace $Y\subseteq X$ with $k=\dim Y\ge c(\e)\log n$ such
that $d(Y,\ell_2^k)\le 1+\e$. Here $d(\cdot,\cdot)$ denotes the
Banach-Mazur distance and $c(\cdot)$ depends only on $\e$. The
first result of this type appeared in \cite{dvoretzky}, and the
logarithmic lower bound on the dimension is due to V.
Milman~\cite{milmandvoretzky}. If in addition to taking subspaces,
we also allow passing to quotients, the dimension $k$ above can be
greatly improved. This is V. Milman's Quotient of Subspace
Theorem~\cite{milmanqs} (commonly referred to as the QS Theorem),
a precise formulation of which reads as follows:

\begin{theorem}[Milman's QS Theorem \cite{milmanqs}]\label{thm:qs} For every
$0<\delta<1$ there is a constant $f(\delta)\in (0,\infty)$ such
that for every $n$-dimensional normed space $X$ there are linear
subspaces $Z\subseteq Y\subseteq X$ with $\dim(Y/Z)=k\ge
(1-\delta)n$ and $d(Y/Z,\ell_2^k)\le f(\delta)$.
\end{theorem}

Over the past two decades, several theorems in the local theory of
Banach spaces were shown to have non-linear analogs. The present
paper, which is a continuation of this theme, is devoted to the
proof of a natural non-linear analog of the QS Theorem, which we
present below.

\medskip

A mapping between two metric spaces $f:M \rightarrow X$, is called
an embedding of $M$ in $X$. The \emph{distortion} of the embedding
is defined as
\[
\mathrm{dist}(f)=\sup_{\substack{x,y\in M\\x\neq
y}}\frac{d_X(f(x),f(y))}{d_M(x,y)}\cdot \sup_{\substack{x,y\in M\\x\neq
y}}\frac{d_M(x,y)}{d_X(f(x),f(y))}.
\]
The least distortion required to embed $M$ in $X$ is denoted by
$c_X(M)$. When $c_X(M)\leq \alpha$ we say that $M$ $\alpha$-embeds
in $X$. If $\M$ is a class of metric spaces then we denote
$c_\M(M)=\inf_{X\in \M}c_X(M)$.

In order to motivate our treatment of the non-linear QS problem,
we first describe a non-linear analog of Dvoretzky's Theorem,
which is based on the following notion: Given a class $\M$ of
metric spaces, we denote by $R_{\M}(\alpha,n)$ the largest integer
$m$ such that any $n$-point metric space has a subspace of size
$m$ that $\alpha$-embeds into $X\in \M$. When $\M=\{\ell_p\}$ we
use the notations $c_p$ and $R_{p}$. The
 parameter $c_2(X)$ is known as the Euclidean distortion of $X$.
 In~\cite{bfm}
Bourgain, Figiel, and Milman study this function, as a non-linear
analog of Dvoretzky's theorem. They prove

\begin{theorem}[Non-Linear Dvoretzky Theorem \cite{bfm}] \label{thm:bfm86}
For any $\alpha>1$ there exists $C(\alpha)>0$ such that $R_2(\alpha,n) \geq
C(\alpha) \log n$. Furthermore, there exists $\alpha_0>1$ such that
$R_2(\alpha_0,n)=O(\log n)$.
\end{theorem}

In~\cite{blmn1} the metric Ramsey problem is studied
comprehensively. In particular, the following
phase transition is proved.

\begin{theorem}[\cite{blmn1}]
\label{thm:phase} The following two assertions hold true:
\begin{enumerate}
\item\label{item:upper} For every $n\in \N$ and $1<\alpha<2$: \(
c(\alpha)\log n\le R_2(\alpha,n)\le 2\log n+C(\alpha), \) where
$c(\alpha),C(\alpha)$ may depend only on $\alpha$.

\item\label{item:lower} For every $\alpha>2$ there is an integer
$n_0$ such that for $n\ge n_0$: \( n^{c'(\alpha)}\le
R_2(\alpha,n)\le n^{C'(\alpha)}, \) where $c'(\alpha),C'(\alpha)$
depend only on $\alpha$ and $0<c'(\alpha)\leq C'(\alpha)<1$.
\end{enumerate}
\end{theorem}

The following result, which deals with the metric Ramsey problem
for large distortion, was also proved in \cite{blmn1}:

\begin{theorem}[\cite{blmn1}]\label{thm:largeramsey} For every $\e>0$, every $n$-point metric space $X$
contains a subset of cardinality at least $n^{1-\e}$ whose
Euclidean distortion is $O\left(\frac{\log (1/\e)}{\e}\right)$.
\end{theorem}

With these results in mind, how should we formulate a non-linear
analog of the QS Theorem? We now present a natural formulation of
the problem, as posed by Vitali Milman.

\medskip

The linear quotient operation starts with a normed space $X$, and
a subspace $Y\subseteq X$, and partitions $X$ into the cosets
$X/Y=\{x+Y\}_{x\in X}$. The metric on $X/Y$ is given by
$d(x+Y,x'+Y)=\inf\{\|a-b\|;\ a\in x+Y,\ b\in x'+Y\}$. This
operation is naturally generalizable to the context of arbitrary
metric spaces as follows: Given a finite metric space $M$,
partition $M$ into pairwise disjoint subsets $U_1,\ldots, U_k$.
Unlike the case of normed spaces, the function
$d_M(U_i,U_j)=\inf\{d_M(u,v);\ u\in U_i, v\in U_j\}$ is not
necessarily a metric on $\U=\{U_1,\ldots, U_k\}$. We therefore
consider the maximal metric on $\U$ majorized by $d_M$, which is
easily seen to be the geodesic metric given by:
$$
d_{\geo}(U_i,U_j)=\inf \left\{\sum_{r=1}^k d_M(V_r,V_{r-1});\
V_0,\ldots, V_k\in \U,\ \ V_0=U_i,\ \ V_k=U_j\right\}.
$$
This operation clearly coincides with the usual quotient
operation, when restricted to the class of normed spaces. When
considering the QS operation, we first pass to a subset of $M$,
and then construct a quotient space as above. We summarize this
discussion in the following definition:

\begin{definition} Let $M$ be a finite metric space. A $Q$ space of $M$ is
a metric space that can be obtained from $M$ by the following
operation: Partition $M$ into $s$ pairwise disjoint subsets
$U_1,\ldots,U_s$ and equip $\U=\{U_1,\ldots,U_s\}$ with the
geodesic metric $d_{\geo}$. Equivalently, consider the weighted
complete graph whose vertices are $\U$ with edge weights:
$w(U_i,U_j)=d_M(U_i,U_j)$. The metric on $\U$ can now be defined
to be the shortest-path metric on this weighted graph. A $Q$ space
of a subset of $M$ will be called a $QS$ space of $M$. Similarly,
a subspace of a $Q$ space of $M$ will be called a $SQ$ space of
$M$.
\end{definition}

The above notion of a quotient of a metric space is due to M.
Gromov (see Section $1.16_+$ in \cite{greenbible}). The
formulation of the non-linear QS problem is as follows: Given
$n\in \mathbb{N}$ and $\alpha\ge 1$, find the largest $s\in
\mathbb{N}$ such that {\em any} $n$-point metric space $M$ has a
QS space of size $s$ that is $\alpha$ embeddable in $\ell_2$. More
generally, we consider the following parameters:

\begin{definition}\label{def:qfunctions} Let $\M$ be a class of metric spaces. For every
$n\in \mathbb{N}$ and $\alpha\ge 1$ we denote by
$\Q_{\M}(\alpha,n)$ (respectively $\QS_{\M}(\alpha,n),
\SQ_\M(\alpha,n)$) the largest integer $m$ such that every
$n$-point metric space has a $Q$ space (respectively $QS$, $SQ$
space) of size $m$ that $\alpha$-embeds into a member of $\M$.
When $\M=\{\ell_p\}$ we use the notations $\Q_p$, $\QS_p$ and
$\SQ_p$.
\end{definition}
We remark that for finite metric spaces, embeddability into
$\ell_p$ and $L_p$ are the same (see e.g.~\cite{dezalau}).
Additionally, by Dvoretzky's theorem~\cite{dvoretzky}, if $X$ is
an infinite dimensional Banach space then for every $n\in
\mathbb{N}$, $\alpha\ge 1$ and $\e>0$, $\Q_{\{X\}}(\alpha,n)\le
\Q_2(\alpha+\e,n)$ (and similarly for $\QS_{\{X\}}, \QS_{\{X\}}$).

\medskip

In the linear setting there is a natural duality between subspaces
and quotients. In particular, one can replace in Dvoretzky's
theorem the word "subspace" by the word "quotient", and the
resulting estimate for the dimension will be identical. Similarly,
the statement of the QS Theorem remains unchanged if we replace
"quotient of subspace" by "subspace of quotient". In the
non-linear setting these simple observations are no longer clear.
In view of Theorem~\ref{thm:phase}, it is natural to ask if the
same is true for $Q$ spaces. Similarly, it is natural to ask if
the $QS$ and $SQ$ functions behave asymptotically the same. In
this paper we present a comprehensive analysis of the functions
$\Q_2$, $QS_2$ and $SQ_2$. It turns out that the answer to the
former question is no, while the answer to the latter question is
yes. On the other hand, as conjectured by Milman, our results show
that just as is the case in the linear setting, once we allow the
additional quotient operation, the size of the Euclidean spaces
obtained increases significantly.

Below is a summary of our results concerning the $QS$ and $SQ$
problems:

\begin{theorem}\label{thm:phaseqs} For every $1<\alpha<2$ there
are constants $0<c(\alpha), C(\alpha)<1$ such that for every $n\in
\N$,
$$
n^{c(\alpha)}\le \QS_2(\alpha,n), \SQ_2(\alpha,n)\le
n^{C(\alpha)}.
$$
On the other hand, for every $\alpha\ge2$ there is an integer
$n_0$ and there are constants $0<c'(\alpha), C'(\alpha)<1$ such
that for every $n\ge n_0$,
$$
c'(\alpha)n\le \QS_2(\alpha,n),\SQ_2(\alpha,n)\le C'(\alpha)n.
$$
\end{theorem}

As mentioned above, the $Q$ problem exhibits a different behavior.
In fact, we have a double phase transition in this case:

\begin{theorem}\label{thm:phaseq} For every $1<\alpha<\sqrt{2}$
there is a constant $C_1(\alpha)$ such that for every $n\in \N$,
$\Q_2(\alpha,n)\le C_1(\alpha)$. For every $\sqrt{2}<\alpha<2$
there are constants $c(\alpha), C(\alpha)$ such that for every
$n\in \N$, $n^{c(\alpha)}\le \Q_2(\alpha,n)\le n^{C(\alpha)}$.
Finally, for every $\alpha\ge2$ there is an integer $n_0$ and
there are constants $0<c'(\alpha), C'(\alpha)<1$ such that for
every $n\ge n_0$, $ c'(\alpha)n\le \Q_2(\alpha,n)\le C'(\alpha)n
$.

\end{theorem}

In other words, for $\alpha>\sqrt{2}$ the asymptotic behavior of the function
$\Q_2$ is the same as the behavior of the functions $\QS_2$ and $\SQ_2$. We
summarize the qualitative behavior of the size of subspaces, quotients,
quotients of subspaces and subspaces of quotients of arbitrary metric spaces
in Table~\ref{tb:euc-qs}. For aesthetic reasons, in this table we write
$\mathcal{S}_2=R_2$ (i.e. $R_2$ is the "subspace" function). The first row
contains results from \cite{blmn1}. We mention here that the behavior of
$\mathcal{S}_2(2,n)$ remains unknown. Furthermore, we do not know the
behavior of the function $\Q_2$ at $\sqrt{2}$. Finally, we mention that
in~\cite{blmn2} it is shown that $R_2(1,n)=3$ for all $n\ge 3$. We did not
study the functions $\Q_2,\QS_2,\SQ_2$ in the isometric case.

\begin{table}[ht]
\begin{center}
\begin{tabular}{|l||c|c|c|} \hline
  & \multicolumn{3}{c}{Distortion} \vline \\
                 & $(1,\sqrt{2})$  &
                $(\sqrt{2},2)$  & $(2,\infty)$ \\ \hline \hline
$\mathcal{S}_2$
  & \multicolumn{2}{c}{logarithmic} \vline  & polynomial \\ \hline
 $\mathcal{Q}_2$& {constant} &  polynomial &
{proportional}  \\ \hline
 $\mathcal{QS}_2$  &\multicolumn{2}{c}{polynomial} \vline &
 proportional \\ \hline
 $\mathcal{SQ}_2$  &\multicolumn{2}{c}{polynomial} \vline &
 {proportional}  \\ \hline
\end{tabular}
\end{center}
\caption{The qualitative behavior of the Euclidean quotient/
subspace functions, for different distortions.} \label{tb:euc-qs}
\end{table}

\medskip

For large distortions we prove the following analog of
Theorem~\ref{thm:largeramsey}:

\begin{theorem}\label{thm:largeq} For any $n\in \mathbb{N}$ and
$\e>0$, every $n$-point metric space has a $Q$ space of size
$(1-\e)n$ whose Euclidean distortion is $O\big(\log(2/\e)\big)$.
On the other hand, there are arbitrarily large $n$-point metric
spaces every $QS$ or $SQ$ space of which, of size at least
$(1-\e)n$, has Euclidean distortion $\Omega\big(\log(2/\e)\big)$.
\end{theorem}

This result should be viewed in comparison to Bourgain's embedding
theorem~\cite{bourgainembedding}, which states that for every
$n$-point metric spaces $X$, $c_2(X)=O(\log n)$.
Theorem~\ref{thm:largeq} states that if one is allowed to identify
an arbitrarily small proportion of the elements of $X$, it
possible to arrive at a metric space whose Euclidean distortion is
bounded independently of $n$. In fact, Theorem~\ref{thm:largeq} is
proved via a modification of Bourgain's original proof. This is
unlike the situation for the non-linear Dvoretzky problem, since
in~\cite{blmn5} an example is constructed which shows that
Bourgain's embedding method cannot yield results such as
Theorem~\ref{thm:largeramsey}.

\smallskip

Except for a loss in the dependence on $\e$, it is possible to give a more
refined description of the $Q$ spaces obtained in Theorem \ref{thm:largeq}.
Using a different embedding method, we can actually ensure that for every
$\e>0$, every $n$ point metric space has a $Q$ space of size $(1-\e)n$ which
well embeds into an {\em ultrametric}. This is of interest since such spaces
have a simple hierarchically clustered structure, which is best described
through their representation as a {\em hierarchically well-separated tree}
(see Section \ref{section:lowerlarge} for the definition). This special
structure is useful in several algorithmic contexts, which will be discussed
in a forthcoming (Computer Science oriented) paper.

\begin{theorem}\label{thm:qhst} For every $\e>0$ and $n\in
\mathbb{N}$, any $n$ point metric space $X$ contains a subset
$A\subseteq X$ of size at most $\e n$ such that the quotient of
$X$ induced by the partition $\big\{\{a\}\big\}_{a\in X\setminus
A}\cup \{A\}$ is $O\left[\frac{\log (1/\e)}{\e}\right]$ equivalent
to an ultrametric. On the other hand, there are arbitrarily large
$n$-point metric spaces every $QS$ or $SQ$ space of which, of size
at least $(1-\e)n$, cannot be embedded into an ultrametric with
distortion $O(1/\e)$.
\end{theorem}

\medskip

In Section~\ref{section:cube} we study the QS problem for the
hypercube $\Omega_d=\{0,1\}^d$ (although the embedding results
used there may also be of independent interest). The cube-analog
of Theorem~\ref{thm:largeramsey} was studied in \cite{blmn1},
where it was shown that if $B\subset \Omega_d$ satisfies
$c_2(B)\le \alpha$ then $|B|\le C2^{(1-c/\alpha^2)d}$, and on the
other hand there is a subset $B_0\subset \Omega_d$ with Euclidean
distortion at most $\alpha$ and which contains at least
$2^{(1-[\log (c'\alpha)]/\alpha^2)d}$ points (here $c,c',C$ are
positive universal constant). In Section~\ref{section:cube} we
prove the following $QS$ counterpart of this result:

\begin{theorem}\label{thm:cube} There is an absolute constant
$c>0$ such that for all $d\in \N$ and $0<\e<1/2$, every $QS$ space
of $\Omega_d$ containing more than $(1-\e)2^d$ points has
Euclidean distortion at least:
$$
c\sqrt{\frac{\log(1/\e)}{1+\log\left(\frac{d}{\log(1/\e)}\right)}}.
$$
On the other hand, there are $QS$ spaces of $\Omega_d$, of size
greater than $(1-\e)2^d$ whose Euclidean distortion matches this
bound.
\end{theorem}

\medskip

In Section~\ref{section:lipq} we briefly study another notion of
quotient introduced by Bates, Johnson, Lindenstrauss, Preiss and
Schechtman~\cite{bates}, which has been the focus of considerable
attention in the last few years. It turns that this notion of
quotient, while being useful in many contexts, does not yield a
satisfactory non-linear version of the QS theorem (at least for
distortion greater than $2$). Namely, we show that using this
notion of quotient we cannot expect to obtain quotients of
subspaces which are asymptotically larger than what is obtained by
just passing to subspaces (i.e. what is ensured by
Theorem~\ref{thm:phase}).

In order to describe this notion we recall the following standard
notation which will be used throughout this paper. Given a metric
space $M$, $x\in M$ and $\rho>0$, denote $B_M(x,\rho)=\{y\in M;\
d_M(x,y)\le \rho\}$ and $B_M^\circ(x,\rho)=\{y\in M;\ d_M(x,y)<
\rho\}$.

Let $(X,d_X)$ and $(Y,d_Y)$ be metric spaces and $c>0$. A function
$f:X\to Y$ is called $c$-co-Lipschitz if for every $x\in X$ and
every $r>0$, $f\big(B_X(x,r)\big)\supseteq B_Y\big(f(x),r/c\big)$.
The function $f$ is called co-Lipschitz if it is $c$-co-Lipschitz
for some $c>0$. The smallest such $c$ is denoted by $\coLip(f)$. A
surjection $f:X\to Y$ is called a Lipschitz quotient if it is both
Lipschitz and co-Lipschitz. The notion of co-Lipschitz mappings
was introduced by Gromov (see Section 1.25 in~\cite{greenbible}),
and the definition of Lipschitz quotients is due to Bates,
Johnson, Lindenstrauss, Preiss and Schechtman~\cite{bates}. The
basic motivation is the fact that the Open Mapping Theorem ensures
that surjective continuous linear operators between Banach spaces
are automatically co-Lipschitz.

In the context of finite metric spaces these notions only make
sense with additional quantitative control of the parameters
involved. Given $\alpha>0$ and two metric spaces $(X,d_X)$,
$(Y,d_Y)$ we say that $X$ has an $\alpha$-Lipschitz quotient in
$Y$ if there is a subset $Z\subset Y$ and a Lipschitz quotient
$f:X\to Z$ such that $\Lip(f)\cdot \coLip(f)\le \alpha$. The
following definition is the analog of
Definition~\ref{def:qfunctions} in the context of Lipschitz
quotients.

\begin{definition}\label{def:lipqfunction} Let $\M$ be a class of metric spaces. For every
$n\in \N$ and $\alpha\ge 1$ we denote by $\QS_\M^{\Lip}(\alpha,n)$
the largest integer $m$ such that every $n$-point metric space has
a subspace which has an $\alpha$-Lipschitz quotient in a member of
$\M$. When $\M=\{\ell_p\}$ then we use the notation $\QS_p^{\Lip}$
\end{definition}

The main result of Section~\ref{section:lipq} is:

\begin{theorem}\label{thm:lipqs} The following two assertions hold
true:
\begin{enumerate}
\item For every $\alpha>2$ there is an integer $n_0$ such that for
$n\ge n_0$: \( n^{c(\alpha)}\le \QS_2^{\Lip}(\alpha,n)\le
n^{C(\alpha)}, \) where $c(\alpha),C(\alpha)$ depend only on
$\alpha$ and $0<c(\alpha)\leq C(\alpha)<1$. \item For every $1\le
\alpha <2$ there is an integer $n_0$ such that for $n\ge n_0$:
$$
e^{c'(\alpha)\sqrt{\log n}}\le \QS_2^{\Lip}(\alpha,n)\le
e^{C'(\alpha)\sqrt{(\log n)(\log \log n)}},
$$
where $c'(\alpha),C'(\alpha)$ depend only on $\alpha$.
\end{enumerate}
\end{theorem}

Thus, the additional Lipschitz quotient operation only yields an
improvement for distortion smaller than $2$. We have not studied
the analogous questions for the $Q$ and $SQ$ problems.

\bigskip

Throughout this paper we also study the functions $\Q_p,\SQ_p,\QS_p$ for
general $1\le p<\infty$. In most cases we obtain matching or nearly matching
upper and lower bounds for the various functions, but some interesting
problems remain open. We summarize in Table~\ref{tb:p<2-qs} and
Table~\ref{tb:p>2-qs} the qualitative nature of our results (in which we
write once more $R_p=\mathcal{S}_p$). As is to be expected, it turns out that
there is a difference between the cases $1\le p\le 2$ and $p>2$. In both
tables, the first row contains results from \cite{blmn1} and \cite{blmn2}. In
Table~\ref{tb:p>2-qs} the question marks refer to the fact that for $p>2$ our
lower and upper bounds do not match in the range $(2^{2/p},2)$.

\begin{table}[ht]
\begin{center}
\begin{tabular}{|l||c|c|c|} \hline
  & \multicolumn{3}{c}{Distortion} \vline \\
 &  $(1,2^{1-\frac{1}{p}})$ &  $(2^{1-\frac{1}{p}},2)$  & $(2,\infty)$ \\ \hline \hline
$\mathcal{S}_p$
   &  \multicolumn{2}{c}{logarithmic} \vline   & polynomial \\ \hline
 $\mathcal{Q}_p$& {constant}  & polynomial &
{proportional}  \\ \hline
 $\mathcal{QS}_p$ & \multicolumn{2}{c}{polynomial} \vline &
 {proportional}  \\ \hline
 $\mathcal{SQ}_p$ & \multicolumn{2}{c}{polynomial} \vline &
 {proportional} \\ \hline
\end{tabular}
\end{center}
\caption{The qualitative behavior of the  $\ell_p$ quotient/
subspace function for $p\leq 2$, and different distortions.}
\label{tb:p<2-qs}
\end{table}

\begin{table}[ht]
\begin{center}
\begin{tabular}{|l||c|c|c|c|} \hline
  & \multicolumn{4}{c}{Distortion} \vline \\
 &  $(1,2^{\frac{1}{p}})$ &  $(2^{\frac{1}{p}},2^{\frac{2}{p}})$
 & $(2^{\frac{2}{p}},2)$ & $(2,\infty)$ \\ \hline \hline
$\mathcal{S}_p$
  &  \multicolumn{2}{c}{logarithmic} \vline & ? & polynomial \\ \hline
 $\mathcal{Q}_p$& {constant}  & {polynomial}  & {?} &  {proportional}  \\ \hline
 $\mathcal{QS}_p$ & \multicolumn{2}{c}{polynomial} \vline & {?} &  {proportional}  \\ \hline
 $\mathcal{SQ}_p$ & \multicolumn{2}{c}{polynomial} \vline & ? & {proportional} \\ \hline
\end{tabular}
\end{center}
\caption{The qualitative behavior of the  $\ell_p$ quotient/
subspace function for $p\geq 2$, and different distortions.}
\label{tb:p>2-qs}
\end{table}

\bigskip

This paper is organized as follows. Section~\ref{section:upper}
deals with the various upper bounds for $\Q_p,\QS_p, \SQ_p$. In
Section \ref{section:lowerlarge} we prove Theorem~\ref{thm:largeq}
and Theorem~\ref{thm:qhst}. In Section~\ref{section:lowersmall} we
prove the various lower bounds for $\Q_p(\alpha,n)$ ,
$\QS_p(\alpha,n)$ and $\SQ_p(\alpha,n)$ for $\alpha\le 2$.
Section~\ref{section:cube} deals with the $QS$ problem for the
hypercube (in the context of embeddings into $\ell_p$ for general
$p\ge 1$). Finally, Section~\ref{section:lipq} deals with the $QS$
problem for Lipschitz quotients.

\section{Upper Bounds}\label{section:upper}

In this section we present the various upper bounds for the $Q$,
$QS$ and $SQ$ problems presented in the introduction. In the
following two sections we will provide matching lower bounds for
these problems.

We begin with an abstract method with which one can obtain upper
bounds for $\QS_\M(\alpha,n)$, for various classes of metric
spaces $\M$.

\begin{lemma}\label{lem:upperexpander} Let $\M$ be a class of
metric spaces and $\alpha>1$. Assume that there exists a $k$-point
metric space $X$ such that $c_\M(X)>\alpha$. Then for every
integer $n$,
$$
\max\left\{\SQ_\M(\alpha,nk),\QS_\M(\alpha,nk)\right\}\le\left(k-\frac12\right)n.
$$
\end{lemma}

\begin{proof} Define
$Y=X\times\{1,\ldots,n\}$. We equip $Y$ with the following metric:
\begin{eqnarray*}
d_Y((x,i),(y,j))=\left\{\begin{array}{ll} d_X(x,y) & i=j\\
\beta & i\neq j
\end{array}\right.
\end{eqnarray*}
It is straightforward to verify that provided $\beta\ge \diam(X)$,
$d_Y$ is indeed a metric.

Since $|Y|=nk$, it is enough to show that $Y$ has no $QS$ or $SQ$
space of size greater than $\left(k-\frac12\right)n$ which
$\alpha$-embeds into a member on $\M$. Let $U_1\ldots,U_r\subseteq
Y$ be disjoint subsets and $r>\left(k-\frac12\right)n$. Denote
$m=|\{1\le i\le r;\ |U_i|=1\}|$. Then:
$$
kn\ge \left|\bigcup_{i=1}^rU_i\right|=\sum_{i=1}^r |U_i|\ge m
+2(r-m)>2\left(k-\frac12\right)n-m=2kn-n-m.
$$
Hence $m>kn-n$, which implies that there is $i\in \{1,\ldots,n\}$ such that
the singletons $\big\{\{(x,i)\}\big\}_{x\in X}$ are all elements of
$\U=\{U_1,\ldots,U_r\}$. If $Y$ has either a $QS$ space or a $SQ$ space of
size greater than $r$ which $\alpha$-embed into a member of $\M$ then we
could find such $\U$ which could be completed to a partition $\V$ of a subset
$S\subseteq Y$ such that $\U$, equipped with the quotient metric induced by
$\V$, $\alpha$-embeds into a member of $\M$. By taking $\beta=\diam(X)$ we
guarantee that both the QS and the SQ metrics induced by $\V$, when
restricted to  $X\times \{i\}$ are isometric to $X$. This contradicts the
fact that $c_\M(X)>\alpha$.
\end{proof}

The next two corollaries are the upper bounds contained in Theorem
\ref{thm:largeq} and the second part of Theorem \ref{thm:phaseqs}.

\begin{corollary}\label{coro:upperbourgain} For every $\e\in
(0,1)$ and $1\le p<\infty$ there are arbitrarily large $n$-point
metric spaces every $QS$ or $SQ$ space of which, $\U$, of size at
least $(1-\e)n$, satisfies $c_p(\U)\ge
\Omega\big([\log(2/\e)]/p\big)$.
\end{corollary}

\begin{proof} By \cite{matexpander} there are constants $c,\e_0>0$ such that for $\e\le \e_0$
there is a $k$-point metric spaces $X$ with $ k\le \frac{1}{3\e}$,
for which $c_p(X)\ge c[\log(1/\e)]/p $. By Lemma
\ref{lem:upperexpander}, for every integer $m$ there is a metric
space of size $km$, such that every $QS$ or $SQ$ space of which,
of size at least $\left(k-\frac{1}{3}\right)m\le (1-\e)km$, cannot
be embedded in $\ell_p$ with distortion smaller than
$c[\log(2/\e)]/p$.
\end{proof}

\begin{corollary}\label{coro:largephase} For every $\alpha>1$ there exists a constant
$c(\alpha)<1$ such that for every $1\le p<\infty$ there is an
integer $n_0=n_0(p)$ such that for every $n\ge n_0$,
$\QS_p(\alpha,n),\SQ_p(\alpha,n)\le c(\alpha)\cdot n$.
\end{corollary}

\begin{proof} By \cite{matexpander} there is a constant $c>0$ such
that for every $k$ large enough there is a metric space $X_k$ such
that for every $1\le p<\infty$, $c_p(X_k)\ge [c\log k]/p$. So, for
$k=\left\lfloor e^{\alpha p/c}\right\rfloor+1$, $c_p(X_k)>\alpha$.
If $n>8k^2$ then we can find an integer $m$ such that
$\frac{n}{k}\le m\le \frac{4k-1}{4k-2}\cdot\frac{n}{k}$. By
Lemma~\ref{lem:upperexpander},
$$
\max\left\{\QS_p(\alpha, n),\SQ_p(\alpha,n)\right\}\le
\max\left\{\QS_p(\alpha, mk),\SQ_p(\alpha,mk)\right\}\le
\left(k-\frac{1}{2}\right)m\le\left(1-\frac{1}{4k}\right)n.
$$
\end{proof}

The upper bound for embedding into the class of ultrametrics, analogous to
Corollary~\ref{coro:upperbourgain}, shows that in this case the asymptotic
dependence on $\e$ is worse. In order to prove it we need the following
simple lemma. Recall that a metric space $(X,d)$ is called ultrametric if for
every $x,y,z\in X$, $d(x,y)\le \max\{d(x,z),d(y,z)\}$. In what follows we
denote by $\mathrm{UM}$ the class of all ultrametrics.

\begin{lemma}\label{lem:lineUM} Let $\{a_i\}_{i=1}^n$ be an
increasing sequence of real numbers, equipped with the metric
induced by the real line. Then:
$$
c_{\mathrm{UM}}(\{a_1,\ldots,a_n\})\ge \frac{a_n-a_1}{\max_{1\le
i\le n-1}(a_{i+1}-a_i)}.
$$
In particular, $c_{\mathrm{UM}}(\{1,\ldots,n\})\ge n-1$, i.e. the
least distortion embedding of $\{1,\ldots,n\}$ into an ultrametric
is an embedding into an equilateral space.
\end{lemma}

\begin{proof} Let $X$ be an ultrametric and
$f:\{a_1,\ldots,a_n\}\to X$ be an embedding such that for all
$1\le i,j\le n$, $d_X(f(a_i),f(a_j))\ge |a_i-a_j|$ and there exist
$1\le i<j\le n$ for which $d_X(f(a_i),f(a_j))= |a_i-a_j|$. For
$1\le i,j\le n$ write $i\sim j$ if $d_X(f(a_i),f(a_j))<a_n-a_1$.
The fact that $X$ is an ultrametric implies that $\sim$ is an
equivalence relation. Moreover, our assumption of $f$ implies that
$1\not\sim n$. It follows that there exists $1\le i\le n-1$ such
that $a_i\not\sim a_{i+1}$, i.e. $d_X(f(a_{i+1}),f(a_i))\ge
a_n-a_1$, which implies the lower bound on the distortion of $f$.
\end{proof}

\begin{corollary}\label{coro:upperUM} For every $0<\e<1$ there are
arbitrarily large $n$-point metric spaces every $QS$ or $SQ$ space
of which, $\U$, of size at least $(1-\e)n$, satisfies
$c_{\mathrm{UM}}(\U)\ge \left\lfloor\frac{1}{2\e}\right\rfloor-2$.
Additionally, for every $\alpha\ge 1$ and every $n\ge
8(\lfloor\alpha\rfloor+2)^2$, $\SQ_{\mathrm{UM}}(\alpha,n),
\QS_{\mathrm{UM}}(\alpha,n)\le \frac{7+4\lfloor
\alpha\rfloor}{8+4\lfloor \alpha\rfloor}n$.
\end{corollary}

\begin{proof} The proof is analogous to the proofs of
Corollary~\ref{coro:upperbourgain} and
Corollary~\ref{coro:largephase}. In the first case we set
$k=\left\lfloor\frac{1}{2\e}\right\rfloor$ and take $X=\{1,\ldots,
k\}$. By Lemma~\ref{lem:lineUM}, $c_{\mathrm{UM}}(X)\ge k-1>k-2$,
and the required result follows from
Lemma~\ref{lem:upperexpander}. In the second case we set
$k=\lfloor\alpha\rfloor+2$ so that $c_{\mathrm{UM}}(\{1,\ldots,
k\})>\alpha$. We conclude exactly as in the proof of
Corollary~\ref{coro:largephase}.
\end{proof}

\medskip

The following proposition bounds from above the functions $\QS_p$
and $\SQ_p$ for distortions smaller than $2^{\min\{1,2/p\}}$. Our
proof is a modification of the technique used in \cite{blmn2}.

\begin{proposition}\label{prop:upperrandom}
There is an absolute constant $c>0$ such that for every $\delta\in
(0,1)$, and every $n\in N$, if $1\le p\le 2$ then:
$$
\max\left\{\QS_p(2-\delta,n), \SQ_p(2-\delta,n)\right\}\le
n^{1-c\delta^2},
$$
and if $2<p<\infty$ then:
$$
\max\left\{\QS_p(2^{2/p}-\delta,n),
\SQ_p(2^{2/p}-\delta,n)\right\}\le n^{1-cp^2\delta^2},
$$
\end{proposition}

\begin{proof} Fix an integer $m$, and denote by $K_{m,m}$ the
complete bipartite $m\times m$ graph. It is shown in \cite{blmn2}
that:
$$
c_p(K_{m,m}) \ge \left\{ \begin{array}{ll}
2\left(\frac{m-1}{m}\right)^{1/p} & 1\le p\le 2\\
2^{2/p}\left(\frac{m-1}{m}\right)^{1/p}&
2<p<\infty.\end{array}\right.
$$
It follows in particular that for $m=\left\lfloor
\frac{4}{p\delta}\right\rfloor$,
$c_p(K_{m,m})>2^{\min\{1,2/p\}}-\delta$.

Fix $0<q<1$, the exact value of which will be specified later. Let
$G=(V,E)$ be a random graph from $G(n,q)$ (i.e. a graph on $n$
vertices, such that each pair of vertices forms an edge
independently with probability $q$).  Define a metric on $V$ by
requiring that for $u,v\in V$, $u\neq v$, $d(u,v)=1$ if $[u,v]\in
E$ and $d(u,v)=2$ if $[u,v]\notin E$. Fix an integer $s$. Consider
a set of $s$ disjoint subsets of $V$, $\U=\{U_1,\ldots,U_s\}$. We
observe that when $\U$ is viewed as either a  QS or SQ space of
$(V,d)$, in both cases the metrics induced by $\U$ are actually
the same (and equal $\min \{d(x,y);\; x\in U_i,\; y\in U_j\}$).
Denote $\W=\left\{U_i;\ 1\le i\le s,\ |U_i|\le
\frac{2n}{s}\right\}$. Clearly $|\W|\ge s/2$. Without loss of
generality, $\W\supseteq \{U_1,U_2,\ldots,U_{\lceil s/2\rceil}\}$.

For $1\le i<j\le \lceil s/2\rceil$ denote by $\gamma_{ij}$ the
probability that there is an edge between $U_i$ and $U_j$. Clearly
$\gamma_{ij}=1-(1-q)^{|U_i|\cdot|U_j|}$, so that:
$$
q\le \gamma_{ij}\le 1-(1-q)^{(2n/s)^2}.
$$
Since $K_{m,m}$ has $m^2$ edges, the probability that the metric
induced by $\U$ (in both of the $SQ$ and $QS$  cases) on a given
$2m$-tuple in $\{U_1,U_2,\ldots,U_{\lceil s/2\rceil}\}$ coincides
with the metric on $K_{m,m}$ is therefore at least:
$$
q^{m^2}\big[ (1-q)^{(2n/s)^2}\big]^{\binom{2m}{2}-m^2}\ge \big[
q(1-q)^{(2n/s)^2}\big]^{m^2}.
$$
As shown in \cite{blmn2}, there are $\left(\frac{s}{4m}\right)^2$
$2m$-tuples of elements of $\{U_1,U_2,\ldots,U_{\lceil
s/2\rceil}\}$, such that any two intersect in at most one point.
Therefore, the probability that $\U$ does not contain a subspace
isometric to $K_{m,m}$ is at most:
$$
\left\{1-\big[
q(1-q)^{(2n/s)^2}\big]^{m^2}\right\}^{\left(\frac{s}{4m}\right)^2}.
$$
Observe that the number of partitions of $V$ into at least $s$
subsets is $s^n+(s+1)^n+\ldots+n^n\le (n+1)^n$, so that the
probability that all the $s$-point $QS$ (or $SQ$) spaces of
$(V,d)$ contain an isometric copy of $K_{m,m}$, and hence cannot
be embedded into $\ell_p$ with distortion smaller that
$2^{\min\{1,2/p\}}-\delta$, is at least:
$$
1-(n+1)^n\left\{1-\big[
q(1-q)^{(2n/s)^2}\big]^{m^2}\right\}^{\left(\frac{s}{4m}\right)^2}.
$$
We will therefore conclude the proof once we verify that for
$s\approx n^{1-cp^2\delta^2}$, we can choose $q$ such that this
probability is positive. Write $s=n^{1-\eta}$ and
$q=p^2\delta^2n^{-2\eta}$. Then, since $m\le \frac{4}{p\delta}$,
there is an absolute constant $C>0$ such that:
$$
\big[ q(1-q)^{(2n/s)^2}\big]^{m^2}\ge Cn^{-32\eta/(p^2\delta^2)}.
$$
Hence:
\begin{eqnarray*}
1-(n+1)^n\left\{1-\big[
q(1-q)^{(2n/s)^2}\big]^{m^2}\right\}^{\left(\frac{s}{4m}\right)^2}&\ge&
1-(n+1)^n\left[1-Cn^{-32\eta/(p^2\delta^2)}\right]^{p^2\delta^2n^{2-2\eta}/16^2}\\
&\ge& 1-e^{n\log (n+1)-C'p^2\delta^2n^{2-64\eta/(p^2\delta^2)}}>0
\end{eqnarray*}
where we have assumed that $C'p^2\delta^2<1$ (which we are clearly
allowed to do), and chosen $\eta=cp^2\delta^2$ for a small enough
constant $c$.

\end{proof}

We end this section by showing that for $\alpha <
2^{\min\{1/p,1-1/p\}}$ we cannot hope to extract quotients of
metric spaces which embed in $\ell_p$ with distortion $\alpha$ and
that contain more than a bounded number of points. This is quite
easy to see, by considering the {\em star metric} (defined below).
What is perhaps less obvious is that star metrics are the only
obstruction for the existence of unboundedly large quotients of
any sufficiently large metric space, as shown in
Section~\ref{section:lowersmall}.

\medskip

Given an integer $n$  we denote by $\bigstar_n$ the metric on
$\{0,1,\ldots,n\}$ given by $d_{\bigstar_n}(i,0)=1$ for $1\le i\le
n$, and $d_{\bigstar_n}(i,j)=2$ for $1\le i<j\le 2$. The metrics
$\bigstar_n$ are naturally called star metrics.

\begin{lemma}\label{lem:lowerstar} For every integer $n$,
\begin{eqnarray}\label{eq:lowerstar}
c_p(\bigstar_n)\ge \left\{\begin{array}{ll}
2^{1-1/p}\left(1-\frac{1}{n}\right)^{1/p} & 1\le p\le 2\\
2^{1/p}\left(1-\frac{1}{n}\right)^{1/p} & 2\le p< \infty.
\end{array} \right.
\end{eqnarray}
\end{lemma}
\begin{proof} Let $f:\bigstar_n\to \ell_p$ be an embedding such
that for every $x,y\in \bigstar_n$, $$d_{\bigstar_n}(x,y)\le
\|f(x)-f(y)\|_p\le L d_{\bigstar_n}(x,y).$$ We begin with case
$1\le p\le 2$. In this case, as shown in \cite{weston}, for every
$x_1,\ldots, x_n,y_1,\ldots,y_n\in \ell_p$,
$$
\sum_{i=1}^n\sum_{j=1}^n\left(\|x_i-x_j\|_p^p+\|y_i-y_j\|_p^p\right)\le
2\sum_{i=1}^n\sum_{j=1}^n\|x_i-y_j\|_p^p.
$$
Applying this inequality to $x_i=f(i)$ and $y_i=f(0)$ we get that:
$$
n(n-1)2^p\le \sum_{i=1}^n\sum_{j=1}^n\|f(i)-f(j)\|_p^p\le
2\sum_{i=1}^n\sum_{j=1}^n\|f(i)-f(0)\|_p^p\le 2n^2L^p.
$$
This proves the required result for $1\le p\le 2$. For $p\ge 2$ we
apply the same argument, but use the following inequality valid
for every $x_1,\ldots,x_n,y_1,\ldots,y_n\in \ell_p$ (see Corollary
7 in \cite{blmn2}):
$$
\sum_{i=1}^n\sum_{j=1}^n\left(\|x_i-x_j\|_p^p+\|y_i-y_j\|_p^p\right)\le
2^{p-1}\sum_{i=1}^n\sum_{j=1}^n\|x_i-y_j\|_p^p.
$$
\end{proof}

In the following corollary (and also later on in this paper), we
use the convention $\Q_p(\alpha,n)=0$ when $\alpha<1$.

\begin{corollary}\label{coro:upperstar} For every integer $n$ and
every $0<\delta<1$, if $1< p\le 2$ then:
$$
\Q_p\big(2^{1-1/p}(1-\delta)^{1/p},n\big)\le 1+\frac{1}{\delta},
$$
and if $2\le p<\infty$ then:
$$
\Q_p\big(2^{1/p}(1-\delta)^{1/p},n\big)\le 1+\frac{1}{\delta}.
$$
\end{corollary}
\begin{proof} It is straightforward to verify that any $Q$ space
of $\bigstar_{n-1}$ of size $k+1$ is isometric to $\bigstar_k$
(the new "root" will be the class containing the old "root" of the
star). The result now follows from the lower bounds in
Lemma~\ref{lem:lowerstar}.
\end{proof}

\section{Lower Bounds for Large
Distortions}\label{section:lowerlarge}

In this section we study the following problem: Given $\e>0$, what
is the least distortion $\alpha$ such that every $n$ point metric
space has a $Q$ space of size $(1-\e)n$ which $\alpha$ embeds into
$\ell_p$? We prove a lower bound which matches the upper bound
proved in Section~\ref{section:upper}. The proof is based on a
modification of Bourgain's fundamental embedding
method~\cite{bourgainembedding}. Next, we further refine the
structural information on the quotients obtained. Namely, we
construct for arbitrary $n$-point metric spaces quotients of size
$(1-\e)n$ which $c(\e)$-embed into an ultrametric. In fact, in
both cases we obtain the following special kind of quotients:

\begin{definition}\label{def:M/A} Let $M$ be an $n$-point metric space and
$A\subseteq M$. Let $\U$ be the partition of $M$ consisting of $A$
and the elements of $M\setminus A$ as singletons. The $Q$ space of
$M$ induced by $\U$ will be denoted $M/A$. By the definition of
the quotient operation it is easy to verify that for every $x,y\in
M\setminus A$,
\begin{eqnarray}\label{eq:metricU}
d_{M/A}(x,y)=\min\{d_M(x,y),d_M(x,A)+d_M(y,A)\}.
\end{eqnarray}
Additionally, for $x\in M\setminus A$, $d_{M/A}(x,A)=d_M(x,A)$.
\end{definition}

This simple description of the quotients we construct, together
with the fact that we can ensure they have the hierarchical
structure of ultrametrics, has algorithmic significance, which
will be pursued in  a future paper.

The following definition will be useful:

\begin{definition}\label{def:mcenter} Let $X$ be a metric space,
$x\in X$ and $m\ge 1$. We shall say $x$ is an $m$-{\em center} of
$X$ if for every $y\in X$ and every $r>0$, if $|B_X(y,r)|\ge m$
then $x\in B_X(y,r)$.
\end{definition}

\begin{lemma}\label{lem:mcenter} Let $M$ be an $n$-point metric
space and $0<\e<1$. Then there exists a subset $T\subseteq M$ such
that $|T|\le \e n$ and  $T$ is a $\frac{2\log(2/\e)}{\e}$-center
of $M/T$.

\end{lemma}

\begin{proof} Set $m=\frac{2\log(2/\e)}{\e}$. For every $x\in M$ denote
by $\rho_x(m)$ the smallest $\rho>0$ for which $|B_M(x,\rho)|\ge
m$. Choose a random subset $T\subseteq M$ as follows: Let $S$ be
the random subset of $M$ obtained by choosing each point with
probability $\e/2$. Define:
$$
T=S\cup \left\{x\in M;\ S\cap
B_M\big(x,\rho_x(m)\big)=\emptyset\right\}.
$$
Then:
$$
\E|T|=\E|S|+\sum_{x\in M}\Pr\Big[S\cap
B_M\big(x,\rho_x(m)\big)=\emptyset\Big]\le \frac{\e
n}{2}+\left(1-\frac{\e}{2}\right)^mn<\e n.
$$

Denote $\U=M/T$. The proof will be complete once we show that:
\begin{eqnarray*}\label{eq:property}
\forall w\in \U,\  \forall r>0\ \  |B_\U(w,r)|\ge m \
\Longrightarrow \ T\in B_\U(w,r).
\end{eqnarray*}
Indeed, if $w=T$ then there is nothing to prove. Otherwise, assume
for the sake of contradiction that $w=x$ for some $x\in M\setminus
T$ with $d_\U(w,T)=d_M(x,T)>r$. By (\ref{eq:metricU}), for every
$y\in M\setminus T$, $d_\U(x,y)=\min\{d_M(x,y),
d_M(x,T)+d_M(y,T)\}$. In particular, if $d_\U(x,y)\le r$ then
$d_M(x,y)=d_\U(x,y)$. Hence $|B_M(x,r)|\ge m$, so that by the
construction of $T$, $T\cap B_M(x,r)\neq \emptyset$, contrary to
our assumption.
\end{proof}

The following lemma shows that metric spaces with an $m$-center
well embed into $\ell_p$. The proof is essentially a repetition of
Bourgain's original argument \cite{bourgainembedding} (we actually
follow Matou{\v{s}}ek's $\ell_p$- variant of Bourgain's theorem
\cite{matexpander}).

\begin{lemma}\label{lem:bourgainlp} Fix $m\ge 1$ and let $X$ be a
metric space which has an $m$-center. Then for every $1\le
p<\infty$,
$$
c_p(X)\le 96 \left\lceil \frac{\log m}{p}\right\rceil.
$$

\end{lemma}

\begin{proof} Let $x$ be an $m$-center of $X$. Set
$q=\left\lceil \frac{\log m}{p}\right\rceil$. Fix $u,v\in X$,
$u\neq v$. For $i\in \{0,1,\ldots, q\}$ let $r_i$ be the smallest
radius such that $|B_X(u,r_i)|\ge e^{pi}$ and $|B_X(v,r_i)|\ge
e^{pi}$. Observe that by the definition of $q$,
$|B_X(u,r_q)|,|B_X(v,r_q)|\ge e^{pq}\ge m$, so that since $x$ is
an $m$-center of $X$, $x\in B_X(u,r_q)\cap B_X(v,r_q)$. This
implies that $r_q\ge \frac{d_X(u,v)}{2}$. Fix $i\in \{1,\ldots,
q\}$. By the definition of $r_i$ we may assume without loss of
generality that $|B^\circ_X(u,r_i)|\le e^{pi}$. If $A\subseteq X$
is such that $A\cap B^\circ_X(u,r_i)=\emptyset$ and $A\cap
B_X(v,r_{i-1})\neq\emptyset$ then $d_X(u,A)-d_X(v,A)\ge
r_i-r_{i-1}$. If $A$ is chosen randomly such that each point of
$X$ is picked independently with probability $e^{-pi}$ then the
probability of the former event is at least:
$$
\left[1-\left(1-\frac{1}{e^{pi}}\right)^{|B_X(v,r_{i-1})|}\right]\cdot\left(1-\frac{1}{e^{pi}}\right)^{|B_X^\circ(u,r_{i})|}\ge
\left[1-\left(1-\frac{1}{e^{pi}}\right)^{e^{p(i-1)}}\right]\cdot\left(1-\frac{1}{e^{pi}}\right)^{e^{pi}}\ge
\frac{1}{8e^p}.
$$
For $A\subseteq X$, denote by $\pi_i(A)$ the probability that a
random subset of $A$, with points from $X$ picked independently
with probability $e^{-pi}$, equals $A$. The above reasoning
implies that:
$$
\sum_{A\subseteq X} \pi_i(A)|d_X(u,A)-d_X(v,A)|^p\ge
\frac{(r_i-r_{i-1})^p}{8e^p},
$$
so that if we define $\alpha_A=\frac{1}{q}\sum_{i=1}^q\pi_i(A)$
then:
\begin{eqnarray*}
\sum_{A\subseteq X} \alpha_A |d_X(u,A)-d_X(v,A)|^p&\ge&
\frac{1}{8qe^p}\sum_{i=1}^q(r_i-r_{i-1})^p\\
&\ge&
\frac{1}{8q^pe^p}\left(\sum_{i=1}^q(r_i-r_{i-1})\right)^p\\&=&
\frac{(r_q-r_0)^p}{8q^pe^p}\ge
\frac{[d_X(u,v)]^p}{16\cdot2^pq^pe^p}.
\end{eqnarray*}

Now, the embedding of $X$ sends an element $u\in X$ to a vector
indexed by the subsets of $X$, such that the coordinate
corresponding to $A\subseteq X$ is $\alpha_A^{1/p}\cdot d_X(u,A)$.
Since $\sum_{A\subseteq X} \alpha_A=1$, such a mapping is
obviously non-expanding, and the above calculation shows that the
Lipschitz constant of its inverse is at most $16^{1/p}\cdot2eq\le
96 q$, as required.

\end{proof}

The following corollary is a direct consequence of
Lemma~\ref{lem:mcenter} and Lemma~\ref{lem:bourgainlp}:

\begin{corollary}\label{coro:bourgainlp} There is an absolute
constant $c>0$ such that for every $1\le p<\infty$ and every
$0<\e<1$, any $n$-point metric space $M$ has a subset $A\subseteq
M$ such that $|A|<\e n$ and:
$$
c_p(M/A)\le 1+\frac{c}{p}\log\left(\frac{2}{\e}\right).
$$

\end{corollary}

\medskip

We can also apply Lemma~\ref{lem:mcenter} to obtain quotients
which embed into ultrametrics. The basic fact about ultrametrics,
already put to good use in \cite{blmn1}, is that they are
isometric to subsets of Hilbert space. Another useful trait of
\emph{finite} ultrametrics is that they have a natural
representation as {\em hierarchically well-separated trees}
(HSTs). We recall the following useful definition, due to Y.
Bartal \cite{bartal}:

\begin{definition}\label{def:hst} Given $k\ge 1$, a $k$-HST is a metric space whose
elements are leaves of a rooted tree $T$. To each vertex $u\in T$,
a label $\Delta(u)$ is associated such that $\Delta(u)=0$ if and
only if $u$ is a leaf of $T$. The labels are strongly decreasing
in the sense that $\Delta(u)\le \Delta(v)/k$ whenever $u$ is a
child of $v$. The distance between two leaves $x,y\in T$ is
defined as $\Delta(\lca(x,y))$, where $\lca(x,y)$ denotes the the
least common ancestor of $x$ and $y$ in $T$. In what follows, $T$
is called the defining tree of the $k$-HST. For simplicity we call
a $1$-HST a HST. It is an easy fact to verify that the notion of a
finite ultrametric coincides with that of a HST. Although $k$-HSTs
will not appear in this section, this proper subclass of
ultrametrics will play a key role in
Section~\ref{section:lowersmall}.

\end{definition}

\begin{lemma}\label{lem:mcenterhst} Let $m\ge 1$ be an integer and let $X$ be an $n$-point metric
space which has an $m$-center. Then $X$ $2m$-embeds into an
ultrametric.

\end{lemma}

\begin{proof} We prove by induction on $n$ that there is a HST $H$ with $\diam(H)=\diam(X)$ and a bijection
$f:X\to H$ such that for every $u,v\in X$, $d_X(u,v)\le
d_H(f(u),f(v))\le 2md_X(u,v)$. For $n=1$ there is nothing to
prove. Assuming $n>1$, let $x$ be an $m$-center of $X$. Denote
$\Delta=\diam(X)$, and let $a,b\in X$ be such that
$d_X(a,b)=\Delta$. We may assume without loss of generality that
$d_X(x,a)\ge \Delta/2$. For every $k=1,\ldots,m$, define:
$$
A_i=\left\{y\in X;\ \frac{\Delta(i-1)}{2m}\le
d_X(y,a)<\frac{\Delta i}{2m}\right\}.
$$
Now, $\cup_{i=1}^m A_i=B^\circ_X(a,\Delta/2)=\{y\in X;\
d_X(a,y)<\Delta/2\}$. Since $X$ is finite, there is some
$r<\Delta/2$ such that $B^\circ_X(a,\Delta/2)=B_X(a,r)$. But
$x\notin B_X(a,r)$, and since $x$ is an $m$-center of $X$, it
follows that $|B_X(a,r)|<m$. Since the set $\{A_i\}_{i=1}^m$ are
disjoint, and $A_1\neq \emptyset$, it follows that there exists
$1\le i\le m-1$ for which $A_{i+1}=\emptyset$.

Denote $B=\cup_{j=1}^i A_j=B^\circ_X(a, \Delta i/(2m))$. Observe
that $X\setminus B$ has an $m$-center (namely $x$), and $B$ has an
$m$-center vacuously (since $|B|<m$). By the inductive hypothesis
there are HSTs $H_1,H_2$, defined by trees $T_1,T_2$,
respectively,  such that $\diam(H_1)=\diam(B)$,
$\diam(H_2)=\diam(X\setminus B)$, and there are bijections
$f_1:B\to H_1$, $f_2:X\setminus B\to H_2$ which are
non-contracting and $2m$-Lipschitz. Let $r_1,r_2$ be the roots of
$H_1,H_2$, respectively. Let $T$ be the labelled tree $T$ rooted
at $r$ such that $\Delta(r)=\diam(X)=\Delta$, $r_1,r_2$ are the
only children of $r$, and the subtrees rooted at $r_1,r_2$ are
isomorphic to $H_1,H_2$, respectively. Since
$\Delta(r_1)=\diam(H_1)=\diam(B)\le \diam(X)$, and similarly for
$r_2$, $T$ defines a HST on its leaves $H=H_1\cup H_2$. We define
$f:X\to H$ by $f|_{B}=f_1$, and $f|_{X\setminus B}=f_2$. If $u\in B$ and $v\in X\setminus B$
then $d_H(f(u),f(v))=\Delta(r)=\Delta\ge d_X(u,v)$. Furthermore,
$d_X(u,a)<\Delta i/2m$ and $d_X(v,a)\ge \Delta(i+1)/2m$ (since
$A_{i+1}=\emptyset$). Hence:
$$
d_X(u,v)\ge
d_X(v,a)-d_X(u,a)>\frac{\Delta}{2m}=\frac{d_H(f(u),f(v))}{2m}.
$$
This concludes the proof.
\end{proof}

Lemma~\ref{lem:mcenter} and Lemma~\ref{lem:mcenterhst} imply the
following corollary:

\begin{corollary}\label{coro:largeQhst} For every $0<\e<1$ and
every integer $n$, every $n$-point metric space $M$ contains a
subset $A\subseteq M$ such that $|A|\le \e n$ and:
$$
c_{\mathrm{UM}}(M/A)\le \frac{6\log(2/\e)}{\e}.
$$

\end{corollary}

\section{Lower Bounds for Small Distortions}\label{section:lowersmall}

In this section we give lower bounds for $\Q_p(\alpha,n)$, $\QS_p(\alpha,n)$,
$\SQ_p(\alpha,n)$ when $\alpha\le 2$. We begin by showing that for distortion
$\alpha$ greater than $2^{\min\{1-1/p,1/p\}}$, every $n$-point metric space
has a polynomially large $Q$ space which $\alpha$-embeds in $\ell_p$. The
following combinatorial lemma will be used several times in this section. In
what follows, given an integer $n\in \N$ we use the notation
$[n]=\{1,\ldots,n\}$. We also denote by $\binom{[n]}{2}$ the set of all
unordered pairs of distinct integers in $[n]$.

\begin{lemma}\label{lem:coloring} Fix $n,k\in \N$, $n\ge 2$.
For every function $\chi:\binom{[n]}{2}\to [k]$ there is an
integer $s\ge \left\lfloor\frac{n^{1/k}}{8\log n}\right\rfloor$
and there are nonempty disjoint subsets $A_1,\ldots,A_s\subseteq
\{1,\ldots,n\}$ and $\ell\in \{1,\ldots,k\}$ such that for every
$1\le i<j\le s$,
$$
\min\left\{ \chi(p,q);\ p\in A_i,\ q\in A_j\right\}=\ell.
$$
Furthermore, for every $1\le i, j\le s$, $i\neq j$, and every
$p\in A_i$, there exists $q\in A_j$ such that $\chi(p,q)=\ell$.
\end{lemma}

\begin{proof} The proof is by induction on $k$. For $k=1$ there is
nothing to prove. Assume that $k>1$ and denote $m=|\{(i,j);\
\chi(i,j)=1\}|$. Define $s=\left\lfloor\frac{n^{1/k}}{8\log
n}\right\rfloor$. We first deal with the case $m\ge
\frac12n^{1+1/k}$. For each $i\in \{1,\ldots, n\}$ let $B_i=\{j;\
\chi(i,j)=1\}$. Denote $C=\{i;\ |B_i|\ge n^{1/k}/4\}$. Then:
$$
\frac{n^{1+1/k}}{2}\le m=\sum_{i=1}^n |B_i|\le |C|n
+(n-|C|)\frac{n^{1/k}}{4}\le |C|n+\frac{n^{1+1/k}}{4},
$$
i.e. $|C|\ge n^{1/k}/4$.

Consider a random partition of $C$ into $s$ subsets
$A_1,\ldots,A_s$, obtained by assigning to each $i\in C$ an
integer $1\le j\le s$ uniformly and independently. The partition
$A_1,\ldots,A_s$ satisfies the required result with $\ell=1$ if
$A_1,\ldots,A_s$ are nonempty and for every $1\le u\le s$, every
$i\in A_u$ and every $v\neq u$, $B_i\cap A_v\neq \emptyset$. The
probability that this event doesn't occur is at most:
\begin{eqnarray*}
\sum_{u=1}^s\Pr(A_u=\emptyset)+\sum_{u=1}^s\sum_{i\in
C}\sum_{v=1}^s\Pr(i\in A_u,\ B_i\cap A_v=\emptyset)&=&
s\left(1-\frac{1}{s}\right)^{|C|}+\\&\phantom{\le}& \sum_{u=1}^s
\sum_{i\in C}\sum_{v=1}^s\frac{1}{s}\left(1-\frac{1}{s}\right)^{|B_i|}\\
&\le&
(n+1)s\left(1-\frac{1}{s}\right)^{n^{1/k}/4}\\
&\le& \frac{n^{1+1/k}}{4\log n}\exp\left(-\frac{n^{1/k}}{4}\cdot
\frac{8\log
n}{n^{1/k}}\right)\\
&\le& \frac{1}{4\log n}<1,
\end{eqnarray*}
so that the required partition exists with positive probability.

\smallskip
It remains to deal with the case $m<\frac12n^{1+1/k}$. In this
case consider the set $D=\{i;\ |B_i|< n^{1/k}\}$. Then
$\frac12n^{1+1/k}> m\ge n^{1/k}(n-|D|)$, so that $|D|>n/2$.
Consider the graph on $D$ in which $i$ and $j$ are adjacent if and
only if $\chi(i,j)=1$. By the definition of $D$, this graph has
maximal degree less than $n^{1/k}$, so that it has an independent
set $I\subseteq D$ of size at least $|D|/n^{1/k}>
\frac12n^{1-1/k}$ (to see this, color $D$ with $n^{1/k}$ colors
and take the maximal color class). The fact that $I$ is an
independent set means that for $i,j\in I$, $\chi(i,j)>1$, so that
we may apply the inductive hypothesis to $I$ and obtain the
desired partition of size at least $\left\lfloor
\frac{|I|^{1/{(k-1)}}}{8\log|I|}\right\rfloor$. We may assume that
$n^{1/k}\ge 2e$, since otherwise the required result is vacuous.
In this case the lower bound on $|I|$ implies that we are in the
range where the function $x\mapsto x^{1/(k-1)}/\log x$ is
increasing, in which case:
$$
\frac{|I|^{1/{(k-1)}}}{8\log|I|}\ge
\frac{\left(\frac12n^{1-1/k}\right)^{1/(k-1)}}{8(1-1/k)\log n}\ge
\frac{n^{1/k}}{8\log n},
$$
where we have used the inequality $(1-1/k)2^{1/(k-1)}\le 1$.

\end{proof}

The relevance of Lemma~\ref{lem:coloring} to the $QS$ problem is
clear. We record below one simple consequence of it. Recall that
the {\em aspect ratio} of a finite metric space $M$ is defined as:
\begin{eqnarray}\label{eq:defaspect}
\Phi(M)=\frac{\diam(M)}{\min_{x\neq y} d_M(x,y)}.
\end{eqnarray}

\begin{lemma}\label{lem:aspect} Let $M$ be an $n$-point metric
space and $1<\alpha\le 2$. Then there is a $QS$ space of $M$,
$\U$, which is $\alpha$ equivalent to an equilateral metric space
and:
$$
|\U|\ge \left\lfloor \frac{n^{(\log \alpha)/[2\log
\Phi(M)]}}{8\log n}\right\rfloor.
$$
\end{lemma}
\begin{proof}
By normalization we may assume that $\min_{x\neq y} d_M(x,y)=1$.
We may also assume that $\alpha<\Phi(M)$. Write $\Phi=\Phi(M)$ and
set $k=\left\lfloor \frac{\log \Phi}{\log \alpha}\right\rfloor+1$.
For every $x,y\in M$, $x\neq y$ there is a unique integer
$\chi(x,y)\in [k]$ such that $d_M(x,y)\in
[\alpha^{\chi(x,y)-1},\alpha^{\chi(x,y)})$.
Lemma~\ref{lem:coloring} implies that there are disjoint subsets
$U_1,\ldots,U_s\subset M$ and an integer $\ell\in [k]$ such that
$s\ge \left\lfloor \frac{n^{(\log \alpha)/[2\log \Phi(M)]}}{8\log
n}\right\rfloor$ and for every $1\le i<j\le s$, $d_M(U_i,U_j)\in
[\alpha^{\ell-1},\alpha^{\ell})$. Consider the $QS$ space
$\U=\{U_1,\ldots,U_s\}$, and observe that since $\alpha\le 2$, any
minimal geodesic joining $U_i$ and $U_j$ must contain only two
points (namely $U_i$ and $U_j$). This implies that $\U$ is
$\alpha$-equivalent to an equilateral space.

\end{proof}

In what follows we use the following definition:

\begin{definition}\label{def:neighbor}Let $M$ be a finite metric space. For $x\in M$ we
denote by $r_M(x)$ the distance of $x$ to its closest neighbor in
$M$:
$$
r_M(x)=d_M(x,M\setminus \{x\})=\min\{d_M(x,y);\ y\in M, y\neq x\}.
$$

For $0<a<b$ it will be convenient to also introduce the following
notation:
$$
M[a,b)=\left\{x\in M;\ a\le r_M(x)<b\right\}.
$$
\end{definition}

\medskip

For the sake of simplicity, we denote
$\theta(p)=\min\left\{\frac{1}{p},1-\frac{1}{p}\right\}$.

In the following lemma we use the notation introduced in
Definition~\ref{def:M/A} in Section~\ref{section:lowersmall}.

\begin{lemma}\label{lem:T} Let $M$ be an $n$-point metric space.
Then there exist two subsets $S,T\subseteq M$  with the following
properties:

\medskip
\noindent{\bf a)} $S\cap T=\emptyset$.

\medskip
\noindent{\bf b)} $|T|\ge n/4$.

\medskip
\noindent{\bf c)} For every $x\in T$ and every subset $S\subseteq
W\subseteq M\setminus\{x\}$, $d_M(x,W)=r_M(x)$.

\medskip
\noindent{\bf d)} For every $A\supseteq M\setminus T$ and every
$x,y\in M\setminus A$:
$$
d_{M/A}(x,y)=\min\{d_M(x,y),r_M(x)+r_M(y)\}, \quad
d_{M/A}(x,A)=r_M(x).$$
\end{lemma}

\begin{proof} Choose a random subset
$S\subseteq M$ by picking each point independently with
probability $1/2$. Define:
$$
T=\left\{x\in M\setminus S;\ d_M(x,S)=r_M(x)\right\}.
$$
To estimate the expected number of points in $T$, for every $x\in
M$ denote by $N_x\subseteq M$ the set of all points $y\in M$ such
that $r_M(x)=d_M(x,M\setminus\{x\})=d_M(x,y)$. Then $x\in T$ if an
only if $x\notin S$ and $N_x\cap S\neq \emptyset$. These two
events are independent and their probability is at least $1/2$.
Hence $\E |T|\ge n/4$. Parts {\bf a),b)} and {\bf c)} are now
evidently true. Part {\bf d)} follows from part {\bf c)} due to
(\ref{eq:metricU}).

\end{proof}

Given an integer $n$ and $0<\tau\le2$, we denote by
$\bigstar_n^\tau$ the metric on $\{0,1,\ldots,n\}$ given by
$d_{\bigstar_n^\tau}(i,0)=1$ for $1\le i\le n$, and
$d_{\bigstar_n^\tau}(i,j)=\tau$ for $1\le i<j\le n$. The metrics
$\bigstar_n^\tau$ will also be called star metrics (recall that
when $\tau=2$ we have previously used the notation
$\bigstar_n=\bigstar_n^2$).

\begin{lemma}\label{lem:findstar} Let $M$ be an $n$-point metric
space, $0<a<b<2a$ and $b/a\le \alpha\le 2b/a$. Let $S,T\subseteq
M$ be as in Lemma~\ref{lem:T}. Write $m=|T\cap M[a,b)|$. Then
there is some $0<\tau\le \frac{2b}{a\alpha}\le 2$ and a $Q$ space
of $M$, $\U$, which is $\alpha$ -equivalent to
$\bigstar_{|\U|}^\tau$ and:
$$
|\U|\ge \left\lfloor \frac{m^{(\log\alpha)/6}}{8\log
m}\right\rfloor.
$$
\end{lemma}

\begin{proof}Consider the set $N=T\cap M[a,b)$. By definition,
for every $x,y\in N$, $x\neq y$,
$$
\min\{d_M(x,y),r_M(x)+r_M(y)\}\in [a,2b).
$$
Setting $k=\left\lceil \frac{\log(2b/a)}{\log
\alpha}\right\rceil-1$, it follows that there is a unique integer
$\chi(x,y)\in[0,k]$ such that
\begin{eqnarray}\label{eq:onecolor}
\frac{2b}{\alpha^{\chi(x,y)+1}}\le \min\{d_M(x,y),r_M(x)+r_M(y)\}<
\frac{2b}{\alpha^{\chi(x,y)}}.
\end{eqnarray}
Denote $ s=\left\lfloor \frac{m^{1/(k+1)}}{8\log m}\right\rfloor$,
and apply Lemma~\ref{lem:coloring} to get an integer $\ell\in
[0,k]$ and disjoint subsets $A_1,\ldots, A_s\subseteq N$ such that
for every $1\le i<j\le s$:
$$
\min\{\chi(x,y);\ x\in A_i,\ y\in A_j\}=\ell.
$$
Let $\U$ be the $Q$ space of $M$ whose elements are $A_1,\ldots,
A_s$ and $A_0=M\setminus \cup_{i=1}^s A_i$. The metric on $\U$ is
described in the following claim:

\begin{claim}\label{claim:metricU} For every $1\le i\le s$,
$d_\U(A_i,A_0)\in [a,b)$. Furthermore, for every $1\le i< j\le s$,
$$
\frac{2b}{\alpha^{\ell+1}}\le d_\U(A_i,A_j)\le
\frac{2b}{\alpha^\ell}.
$$
\end{claim}
\begin{proof} By part ${\bf c)}$ of Lemma~\ref{lem:T}, for every $1\le i\le s$ and every
$x\in A_i$, $d_M(x,A_0)=r_M(x)\in [a,b)$. Moreover, for every
$0\le i<j\le s$, $d_M(A_i,A_j)\ge a$. Since $2a>b$, this shows
that any geodesic in $\U$ connecting $A_0$ and $A_i$ cannot
contain more than two elements, i.e.
$d_\U(A_i,A_0)=d_M(A_i,A_0)\in [a,b)$. Now, take any $1\le i<j\le
s$. By (\ref{eq:onecolor}), $d_M(A_i,A_j)\ge
\frac{2b}{\alpha^{\ell+1}}$ and:
\begin{eqnarray}\label{eq:upperU}
d_\U(A_i,A_j)&\le&
\min\{d_M(A_i,A_j),d_M(A_i,A_0)+d_M(A_j,A_0)\}\nonumber\\&=&
\min_{x\in A_i,\ y\in A_j}\min\{d_M(x,y),r_M(x)+r_M(y)\}\in
\left[\frac{2b}{\alpha^{\ell+1}},\frac{2b}{\alpha^\ell}\right).
\end{eqnarray}
Consider a geodesic connecting $A_i$ and $A_j$. It is either
$(A_i,A_j)$, $(A_i,A_0,A_j)$ or else it contains either a
consecutive pair $(A_u,A_v)$ for some $1\le u\le v\le s$, $u\neq
v$, or four consecutive pairs $(A_i,A_0), (A_0,A_u),(A_u,A_0),
(A_0,A_v)$ for some $1\le u,v\le s$. In the first three cases we
get that $d_\U(A_i,A_j)\ge \frac{2b}{\alpha^{\ell+1}}$. The fourth
case can be ruled out since in this case the length of the
geodesic is at least $4a>2b\ge \frac{2b}{\alpha^\ell}$, which is a
contradiction to the upper bound in (\ref{eq:upperU}).
\end{proof}

Setting $\tau=\frac{2b}{a\alpha^{\ell+1}}\le 2$, it follows from
Claim~\ref{claim:metricU} that $\U$ is $\alpha$-equivalent to
$\bigstar_s^\tau$.
\end{proof}

The relevance of the metrics $\bigstar_n^\tau$ to the $Q$ problem
is that when $\tau$ is small enough they isometrically embed into
$L_p$:

\begin{lemma}\label{lem:starinLp} For every integer $n$, every $1\le p\le \infty$ and every
$0<\tau\le 2^{1-\theta(p)}$, $\bigstar_n^\tau$ isometrically
embeds into $L_p$.
\end{lemma}

\begin{proof} We begin with the case $1\le p\le 2$. In this case
our assumption implies that there exists $0\le \delta<1$ such that
$\tau=2^{1/p}(1-\delta)^{1/p}$. Our claim follows from the fact
that there are $w_1,\ldots,w_n\in L_p$ such that $\|w_i\|_p=1$ and
for $i\neq j$, $\|w_i-w_j\|_p=\tau$. Indeed, if $\delta=0$ then we
can take these vectors to be the first $n$ standard unit vectors
in $\ell_p$. For $\delta>0$ we take $w_1,\ldots,w_s$ to be i.i.d.
random variables which take the value $\delta^{-1/p}$ with
probability $\delta$ and the value $0$ with probability
$1-\delta$.

The case $p>2$ is slightly different. In this case our assumption
is that $\tau\le 2^{1-1/p}$, so that we may find $0\le \delta\le
1$ such that $\tau=2^{1+1/p}[\delta(1-\delta)]^{1/p}$. We claim
that there are $w_1,\ldots,w_s\in L_p$ such that $\|w_i\|_p=1$ and
for $1\le i<j\le s$, $\|w_i-w_j\|_p=\tau$. Indeed, we can take
$w_1,\ldots,w_s$ to be i.i.d. random variables which take the
value $+1$ with probability $\delta$ and the value $-1$ with
probability $1-\delta$.

\end{proof}

We will require the following definition from \cite{blmn3}.

\begin{definition}\label{def:lacunary} Fix $k\ge 1$. A metric $d$
on $\{1,\ldots, n\}$ is called $k$-lacunary if there is a sequence
$a_1\ge a_2\ge \ldots\ge a_{n-1}\ge 0$ such that $a_{i+1}\le
a_i/k$ and for $1\le i<j\le n$, $d(i,j)=a_i$.
\end{definition}

It is clear that $k$-lacunary spaces are ultrametrics, so that
they embed isometrically in Hilbert space.

\begin{proposition}\label{prop:Qdichotomy} Let $M$ be an $n$-point
metric space, $n\ge 2$. Fix $k\ge 1$, $1<\beta\le 2$ and
$\beta<\alpha<2\beta$. Then $M$ has a $Q$ space, $\U$, which is
$\alpha$-equivalent to either a $k$-lacunary space or a star
metric $\bigstar_{|\U|}^\tau$ for some $0<\tau\le 2/\beta$, and:
$$
|\U|\ge \frac{1}{32\log
n}\left[\frac{n\log(\alpha/\beta)}{\max\{\log
(2k),\log[1/(\beta-1)]\}}\right]^{(\log\alpha)/8}.
$$

\end{proposition}

\begin{proof}Let $T$ be as in Lemma~\ref{lem:T}. For every integer $i\in
\mathbb{Z}$ set:
$$
C_i=\left\{x\in T;\ \left(\frac{\alpha}{\beta}\right)^i\le
r_M(x)<\left(\frac{\alpha}{\beta}\right)^{i+1}\right\}=T\cap
M[(\alpha/\beta)^i,(\alpha/\beta)^{i+1}).
$$
Define $m=\left\lceil \frac{\max\{\log
k,\log[1/(\beta-1)]\}}{\log(\alpha/\beta)}\right\rceil$. For every
$j\in \{0,1,\ldots,m-1\}$ define $D_j=\cup_{i\equiv j
\mod(m)}C_i$. Let $q$ be such that $|D_q|=\max_{i\in \{0,\ldots,
m-1\}}|D_i|$. Then $|D_q|\ge |T|/m\ge n/(4m)+1$.

Set $\ell=|\{r\in \mathbb{Z};\ C_{q+rm}\neq \emptyset\}|$. There
are $r_1>r_2>\cdots>r_\ell$ such that $C_{q+r_im}\neq \emptyset$.
Fix $v_i\in C_{q+r_im}$. Consider the subset
$A=M\setminus\{v_1,\ldots,v_\ell\}\supseteq M\setminus T$. By our
choice of $T$:
\begin{eqnarray*}\label{eq:inci}
d_{M/A}(v_i,A)=r_M(v_i)\in \left[
\left(\frac{\alpha}{\beta}\right)^{q+r_im},\left(\frac{\alpha}{\beta}\right)^{q+r_im+1}\right),
\end{eqnarray*}
and
\begin{eqnarray*}\label{eq:metriccase1}
d_{M/A}(v_i,v_j)=\min\{d_M(v_i,v_j),r_M(v_i)+r_M(v_j)\}.
\end{eqnarray*}
In particular, since $d_M(v_i,v_j)\ge \max \{r_M(v_i),r_M(v_j)\}$,
$$
d_{M/A}(v_i,v_j)\ge
d_M(v_j,S)\ge\left(\frac{\alpha}{\beta}\right)^{q+r_jm}.
$$
Additionally, for $i<j$, since $r_i\ge r_j+1$,
\begin{eqnarray*}
d_{M/A}(v_i,v_j)&\le& r_M(v_i)+r_M(v_j)\\&\le&
\left(\frac{\alpha}{\beta}\right)^{q+r_im+1}+\left(\frac{\alpha}{\beta}\right)^{q+r_jm+1}\\&\le&
\left(\frac{\alpha}{\beta}\right)^{q+r_im}\frac{\alpha}{\beta}\left[1+\left(\frac{\beta}{\alpha}\right)^m\right]\\
&\le& \alpha\left(\frac{\alpha}{\beta}\right)^{q+r_i m},
\end{eqnarray*}
by our choice of $m$. Denote $a_i=(\alpha/\beta)^{q+r_im}$. Then
$a_{i+1}\le (\beta/\alpha)^ma_i\le a_i/k$ and we have shown that
$M/A$ is $\alpha$-equivalent to the $k$-lacunary induced by
$(a_i)$ on $\{1,\ldots,\ell+1\}$.

Let $r$ be such that $|C_{q+rm}|=\max_{i\in \mathbb{Z}}
|C_{q+im}|$. Then $|C_{q+rm}|\ge |D_q|/\ell\ge n/(4m\ell)$. By
Lemma~\ref{lem:findstar}, $M$ has a $Q$ space $\V$ which is
$\alpha$ equivalent to $\bigstar_{|\V|}^\tau$ for some $0<\tau\le
2/\beta$, and:
$$
|\V|\ge \frac{1}{16\log
n}\left(\frac{n}{4m\ell}\right)^{(\log\alpha)/6}.
$$

Summarizing, we have proved the existence of the required $Q$
space of $M$ whose cardinality is at least:
$$
\min_{\ell\ge 1}\max\left\{\ell,\frac{1}{16\log
n}\left(\frac{n}{4m\ell}\right)^{(\log\alpha)/6}\right\},
$$
from which the required result easily follows.

\end{proof}

\begin{corollary}\label{coro:Qsmall} For every $0<\e<1$ and $1<p<\infty$ there exists a constant
$c=c(p,\e)>0$ such that for every integer $n$,
$$
\Q_p\big(2^{\theta(p)}(1+\e),n\big)\ge
cn^{[\theta(p)+\log(1+\e)]/10}.
$$
\end{corollary}

\begin{proof} Apply Proposition~\ref{prop:Qdichotomy} with $k=1$,
$\alpha=2^{\theta(p)}(1+\e)$ and $\beta=2^{\theta(p)}$.
By Lemma~\ref{lem:starinLp}, the resulting $Q$ space is $\alpha$-equivalent
to a subset of $L_p$.
\end{proof}

\begin{corollary}\label{coro:QSsmall} Fix $0<\e<1$. For every
integer $n\ge 2$ and every $1\le p\le \infty$:
$$
\SQ_p(1+\e,n)\ge \frac{n^{\e/12}}{100\log n}.
$$
In fact, for every $k\ge 1$, any $n$-point metric space $M$ has a
$SQ$ space $\U$ which is $1+\e$ embeddable in either a
$k$-lacunary space or an equilateral space, and:
$$
|\U|\ge \frac{1}{100\log n}\left(\frac{n}{\log
(2k)}\right)^{\e/12}.
$$
\end{corollary}

\begin{proof} Apply Proposition~\ref{prop:Qdichotomy} with
$\alpha=1+\e$ and $\beta=\sqrt{1+\e}$. If the resulting $Q$ space
is a star metric then pass to a $SQ$ space by deleting the root of
the star, so that the remaining space is equilateral.

\end{proof}

Before passing to the $QS$ problem, we show that for distortion
$2$ we can obtain proportionally large $Q$ spaces of arbitrary
metric spaces.

\begin{lemma}\label{lem:Q2} For every integer $n$ and every $1\le
p\le \infty$, $ \Q_p(2,n)\ge \frac{n}{4}+1$.

\end{lemma}

\begin{proof} Let $M$ be an $n$-point metric space and let $T$ be
as in Lemma~\ref{lem:T}. Write $|T|=k$ and consider the $Q$ space
$M/A$, where $A=M\setminus T$. We relabel the elements of $M/A$ by
writing $T=\{1,\ldots, k\}$, $A=k+1$, where $r_M(1)\ge
r_M(2)\ldots\ge r_M(k)$. For every $1\le i\le k$,
$d_{M/A}(i,k+1)=r_M(i)$, and for every $1\le i<j\le k$:
$$
d_{M/A}(i,j)=\min\{d_M(i,j),r_M(i)+r_M(j)\}\in [r_M(i),2r_M(i)].
$$
This shows that $M/A$ is $2$ equivalent to the $1$-lacunary space
induced on $\{1,\ldots,k+1\}$ by the sequence
$\{r_M(i)\}_{i=1}^k$.

\end{proof}

\medskip

As we have seen in the proof Corollary~\ref{coro:QSsmall}, the
reason why the $SQ$ problem is "easier" than the $Q$ problem is
that we are allowed to discard the "root" of $Q$ spaces which are
approximately stars. This "easy solution" is not allowed when
dealing with the $QS$ problem. The solution of the $QS$ problem
for distortions less than $2$ is therefore more complicated, and
the remainder of this section is devoted to it.

Our approach to the $QS$ problem builds heavily on the techniques
and results of \cite{blmn1}. Among other things, as in
\cite{blmn1}, we will approach the problem by tackling a more
general {\em weighted} version of the $QS$ problem, which we now
introduce.

\begin{definition} A weighted metric space $(M,d_M,w)$ is a metric
space $(M,d_M)$ with non-negative weights $w:M\to [0,\infty)$.
Given $A\subseteq M$ we denote $w_\infty(A)=\sup_{x\in A} w(x)$.

Given two classes of metric spaces $\M,\A$, and $\alpha\ge 1$ we
denote by $\sigma_\A(\M,\alpha)$ the largest $\sigma\le 1$ such
that any weighted metric space $(M,d_M,w)\in \M$ has a $QS$ space
$\U$ which is $\alpha$-embeddable in a member of $\A$ and
satisfies:
$$
\sum_{A\in \U}w_\infty(A)^\sigma\ge \left(\sum_{x\in M}
w(x)\right)^\sigma.
$$
When $\A$ is the class of all $k$-HSTs we use the notation
$\sigma_k=\sigma_\A$. The case $w\equiv 1$ shows that lower bounds
for $\sigma_k(\M,\alpha)$ also imply lower bounds for the $QS$
problem.

\end{definition}

Having introduced the weighted $QS$ problem, it is natural that we
require a weighted version of Lemma~\ref{lem:coloring}:

\begin{lemma}\label{lem:weighted coloring} Fix $n,k\in \N$, $n\ge
2$, a function $\chi:\binom{[n]}{2}\to [k]$ and a weight function
$w:[n]\to [0,\infty)$. There are disjoint subsets
$A_1,\ldots,A_s\subseteq \{1,\ldots,n\}$ and $\ell\in
\{1,\ldots,k\}$ such that for every $1\le i<j\le s$,
\begin{eqnarray}\label{eq:equicolor}
\min\left\{ \chi(p,q);\ p\in A_i,\ q\in A_j\right\}=\ell,
\end{eqnarray}
and:
$$
\sum_{i=1}^s w_\infty(A_i)^{1/[8k\log(k+1)]}\ge
\left(\sum_{r=1}^nw(r)\right)^{1/[8k\log(k+1)]}.
$$
\end{lemma}

\begin{proof} We use the following fact proved in \cite{bbm}: Let
$x=\{x_i\}_{i=1}^\infty$ be a sequence of non-increasing
non-negative real numbers. Then there exists a sequence
$y=\{y_i\}_{i=1}^\infty$ such that $y_i\le x_i$ for all $i\ge 1$,
$\sum_{i\ge 1} y_i^{1/2}\ge \left(\sum_{i\ge 1}x_i\right)^{1/2}$,
and either $y_i=0$ for all $i>2$ or there exits $w>0$ such that
for all $i\ge 1$, $y_i\in \{w,0\}$. Applying this fact to the
weight function $w:[n]\to [0,\infty)$ we get in the first case
$i,j\in [n]$ such that $\sqrt{w(i)}+\sqrt{w(j)}\ge
\left(\sum_{r=1}^nw(r)\right)^{1/2}$, and we take $A_1=\{i\}$,
$A_2=\{j\}$, $\ell=\chi(i,j)$. In the second case we find $w>0$
and $A\subseteq [n]$, $|A|\ge 3$, such that for $i\in A$ $w(i)\ge
w$ and $|A|\sqrt{w}\ge \left(\sum_{r=1}^nw(r)\right)^{1/2}$. In
this case we may apply Lemma~\ref{lem:coloring} to $A$ and get an
integer $\ell\in [k]$ and disjoint subsets
$A_1,\ldots,A_s\subseteq A$ satisfying (\ref{eq:equicolor}) and
such that $s\ge \left\lfloor \frac{|A|^{1/k}}{8\log
|A|}\right\rfloor$. We can obviously also always ensure that $s\ge
2$. Hence, using the elementary inequality
$\max\left\{\frac{x^{1/k}}{8\log x},2\right\}\ge
x^{1/[4k\log(k+1)]}$, valid for all $x\ge 3$ and $k\ge 1$, we get
that:
$$
\sum_{i=1}^s w_\infty(A_i)^{1/[8k\log(k+1)]}\ge
\big(|A|\sqrt{w}\big)^{1/[4k\log(k+1)]}\ge
\left(\sum_{r=1}^nw(r)\right)^{1/[8k\log(k+1)]}.
$$

\end{proof}

It is useful to introduce the following notation:

\begin{definition}\label{def:M Phi}
For $\Phi\ge 1$ denote by $\M(\Phi)$ the class of all metric
spaces with aspect ratio at most $\Phi$. The class $\M(1)$
consists of all equilateral metric spaces, and is denoted by
$\text{EQ}$.
\end{definition}

We have as a corollary the following weighted version of
Lemma~\ref{lem:aspect}:

\begin{corollary}\label{coro:weightedaspect}
For every $\Phi\ge 2$ and $1\le \alpha\le 2$,
$$
\sigma_{\mathrm{EQ}}(\M(\Phi), \alpha)\ge\frac{\log
\alpha}{16(\log \Phi)\log\left(\frac{2\log \Phi}{\log
\alpha}\right)}.
$$

\end{corollary}

\medskip

We recall below the notion of metric composition, which was used
extensively in \cite{blmn1}.

\begin{definition}[Metric Composition~\cite{blmn1}]\label{def:metric-composition}
Let $M$ be a finite metric space. Suppose that there is a
collection of disjoint finite metric spaces $N_x$ associated with
the elements $x$ of $M$. Let $\mathcal{N} = \{ N_x \}_{x \in M}$.
For $\beta\geq 1/2$, the $\beta$-composition of $M$ and
$\mathcal{N}$, denoted by $C=M_\beta[\mathcal{N}]$, is a metric
space on the disjoint union $\dot \cup_x N_x$. Distances in $C$
are defined as follows. Let $x,y \in M$ and $u \in N_x, v \in
N_y$, then:
\[ d_C(u,v)= \begin{cases} d_{N_x}(u,v) & x=y \\
  \beta \gamma d_M(x,y) & x\neq y .\end{cases}  \]
where $\gamma=\frac{\max_{z \in M} \diam(N_z)}{\min_{x\neq y \in
M} d_M(x,y)}$. It is easily checked that the choice of the factor
$\beta\gamma$ guarantees that $d_C$ is indeed a metric.
\end{definition}

\begin{definition}[Composition Closure~\cite{blmn1}] \label{def:comp}
Given a class $\M$ of finite metric spaces, we consider
$\comp_\beta(\M)$, its closure under $\ge \beta$-compositions.
Namely, this is the smallest class $\C$ of metric spaces that
contains all spaces in $\M$, and satisfies the following
condition: Let $M \in \M$, and associate with every $x \in M$ a
metric space $N_x$ that is isometric to a space in $\C$. Also, let
$\beta'\geq \beta$. Then $M_{\beta'}[\mathcal{N}]$ is also in
$\C$.

\end{definition}

\begin{lemma}\label{lem:composing} Let $\M$ be a class of metric
spaces, $k\ge 1$, $\alpha> 1$ and $\beta\ge \alpha k$. Then:
$$
\sigma_k\big(\mathrm{comp}_\beta(\M),(1+1/\beta)\alpha\big)\ge
\sigma_k(\M,\alpha),
$$
\end{lemma}

\begin{proof} Set $\sigma=\sigma_k(\M,\alpha)$ and take $X\in
\comp_\beta(\M)$. We will prove that for any $w:X\to [0,\infty)$
there exists a $QS$ space $Y$ of $X$ and a $k$-HST $H$ such that
$Y$ is $\alpha$-equivalent to $H$ via a non-contractive
$(1+1/\beta)\alpha$-Lipschitz embedding, and:
$$
\sum_{x\in Y} w(x)^\sigma\ge \left(\sum_{x\in X}
w(x)\right)^\sigma.
$$

The proof is by structural induction on the metric composition. If
$X\in \M$ then this holds by the definition of $\sigma$.
Otherwise, let $M\in \M$ and $\mathcal{N}=\{N_z\}_{z\in
M}\subseteq \comp_\beta(\M)$ be such that
$X=M_\beta[\mathcal{N}]$.

For every $z\in M$ define $w'(z)=\sum_{u\in N_z} w(u)$. By the
definition of $\sigma$ there are disjoint subsets
$U_1,\ldots,U_s\subseteq M$ such that the $QS$ space of $M$,
$\U=\{U_1,\ldots,U_s\}$, is $\alpha$-equivalent to a $k$-HST
$H_M$, defined by the tree $T_M$, via a non-contractive
$\alpha$-Lipschitz embedding, and:
$$
\sum_{i=1}^sw'_\infty(U_i)^\sigma\ge \left(\sum_{z\in M}
w'(z)\right)^\sigma=\left(\sum_{x\in X} w(x)\right)^\sigma.
$$

By induction for each $z\in M$ there are disjoint subsets
$U^z_1,\ldots,U^z_{s(z)}\subseteq N_z$ such that the $QS$ space of
$N_z$, $\U_z=\{U^z_1,\ldots,U_{s(z)}^z\}$ is
$(1+1/\beta)\alpha$-equivalent to a $k$-HST, $H_z$, defined by the
tree $T_z$, via a non-contractive $(1+1/\beta)\alpha$-Lipschitz
embedding, and:
$$
\sum_{i=1}^{s(z)}w_\infty(U_i^z)^\sigma\ge \left(\sum_{u\in N_z}
w(u)\right)^\sigma.
$$
For every $1\le i\le s$ let $z_i\in M$ be such that
$w'(z_i)=w'_\infty(U_i)$. Define $V_1^{z_i},\ldots,
V_{s(z_i)}^{z_i}\subseteq X$ by:
$$
V_1^{z_i}=U_1^{z_i}\bigcup\left(\bigcup_{z\in U_i\setminus
\{z_i\}}N_{z}\right)\qquad \mathrm{and} \qquad V_j^{z_i}=U_j^{z_i}
\ \ \mathrm{for}\ \  j=2,3,\ldots,s(z_i).
$$

Consider the $QS$ space of $X$: $\V=\{V_j^{z_i};\ i=1,\ldots,s\ \
j=1,\ldots,s(z_i)\}$. First of all:
\begin{eqnarray*}
\sum_{A\in \V}
w_\infty(A)^\sigma&=&\sum_{i=1}^s\sum_{j=1}^{s(z_i)}
\left[\max_{x\in V_j^{z_i}} w(x)\right]^\sigma\ge
\sum_{i=1}^s\sum_{j=1}^{s(z_i)} \left[\max_{x\in U_j^{z_i}}
w(x)\right]^\sigma\\
&\ge& \sum_{i=1}^s\left(\sum_{u\in N_{z_i}}w(u)\right)^\sigma =
\sum_{i=1}^s w'(N_{z_i})^\sigma =\sum_{i=1}^s
w'_\infty(U_i)^\sigma \ge \left(\sum_{x\in X} w(x)\right)^\sigma.
\end{eqnarray*}
Therefore, all that remains is to show that $\V$ is
$(1+1/\beta)\alpha$-equivalent to a $k$-HST via a non-contractive,
$(1+1/\beta)\alpha$-Lipschitz embedding. For this purpose we first
describe the metric on $\V$:

\begin{claim}\label{claim:metriccomp} For every $1\le i\le s$ and
every $1\le p< q\le s(z_i)$,
\begin{eqnarray}\label{eq:inblock}
d_\V(V_p^{z_i},V_q^{z_i})=d_{\U_{z_i}}(U_p^{z_i},U_q^{z_i}).
\end{eqnarray}
Furthermore, for every $1\le i<j\le s$ and every $1\le p\le
s(z_i)$, $1\le q\le s(z_j)$:
\begin{eqnarray}\label{eq:betweenblocks}
\beta\gamma d_\U(U_i,U_j)\le d_\V(V_p^{z_i},V_q^{z_j})\le
(\beta+1)\gamma d_\U(U_i,U_j).
\end{eqnarray}
\end{claim}

\begin{proof} By the definition of metric composition, if $z\neq
z_i$, $u\in N_{z_i}$, $v\in N_z$, then $d_X(u,v)\ge
\beta\diam(N_{z_i})>\diam(N(z_i))$. Since for every $1\le j\le
s(z_i)$, $N_{z_i}\cap V_j^{z_i}=U_j^{z_i}$, this implies that
$d_X(V_p^{z_i},V_q^{z_i})=d_{N_{z_i}}(U_p^{z_i},U_q^{z_i})$. In
particular, it follows that $d_\V(V_p^{z_i},V_q^{z_i})\le
d_{\U_{z_i}}(U_p^{z_i},U_q^{z_i})\le \diam(N_{z_i})$. A geodesic
connecting $V_p^{z_i}$ and $V_q^{z_i}$ in $\V$ cannot go out of
$\{V_1^{z_i},\ldots, V_{s(z_i)}^{z_i}\}$, since by the above
observation it would contain a step of length greater that
$\diam(N_{z_i})$. This concludes the proof of (\ref{eq:inblock}).

Next take $1\le i,j\le s$ and $1\le p\le s(z_i)$, $1\le q\le
s(z_j)$ and observe that $d_X(V_p^{z_i},V_q^{z_j})\ge \beta \gamma
d_M(U_i,U_j)$. Indeed, if $i=j$ there is nothing to prove, and if
$i\neq j$ then this follows from the definition of metric
composition and the fact that $V_p^{z_i}\subseteq \cup_{z\in
U_i}N_z$ and $V_p^{z_j}\subseteq \cup_{z\in U_j}N_z$. This
observation implies the left-hand side inequality in
(\ref{eq:betweenblocks}).

To prove the right-hand side inequality in
(\ref{eq:betweenblocks}), take a geodesic
$U_i=W_0,W_1,\ldots,W_m=U_j\in \U$ such that $m$ is minimal. This
implies that $W_r\neq W_{r-1}$ for all $r$, and:
$$
d_\U(U_i,U_j)=\sum_{r=1}^md_M(W_{r-1},W_r),
$$
Let $a_r\in W_{r-1}$, $b_r\in W_r$ be such that
$d_M(a_r,b_r)=d_M(W_{r-1},W_r)$. By construction, for each $r$
there are $A_r,B_r\in \V$ such that $A_r\subseteq N_{a_r}$ and
$B_r\subseteq N_{b_r}$. Consider the following path in $\V$
connecting $V_p^{z_i}$ and $V_q^{z_j}$:
$\Gamma=(V_p^{z_i},A_1,B_1,A_2,B_2,\ldots,
A_{m},B_{m},V_q^{z_j})$. Observe that since $V_p^{z_i}, A_1$
contain points from $N_{z_i}$ and $A_1,B_1$ do not contain points
from a common $N_z$, the definition of metric composition implies
that $d_X(A_1,B_1)\ge \beta d_X(V_p^{z_i},A_1)$. In other words,
$d_X(V_p^{z_i},A_1)+d_X(A_1,B_1)\le
(1+1/\beta)d_X(A_1,B_1)=(\beta+1)\gamma d_M(W_{0},W_1)$.
Similarly, for $r\ge 2$, $d_X(B_{r-1},A_r)+d_X(A_r,B_r)\le
(\beta+1)\gamma d_M(W_{r-1},W_r)$ and
$d_X(A_m,B_m)+d_X(B_m,V_q^{z_j})\le (\beta+1)\gamma d_M(A_m,B_m)$.
Hence, the length of $\Gamma$ is at most $(\beta+1)\gamma
d_\U(U_i,U_j)$, as required.

\end{proof}

We now construct $H$ a $k$-HST that is defined by a tree $T$, as
follows. Start with a tree $T'$ that is isomorphic to $T_M$ and
has labels $\Delta(u) = (\beta+1)\gamma \cdot \Delta_{T_M}(u)$. At
each leaf of the tree corresponding to a point $U_i \in \U$,
create a labelled subtree rooted at $U_i$ that is isomorphic to
$T_{z_i}$ with labels as in $T_{z_i}$. Denote the resulting tree
by $T$. Since we have a non-contractive
$(1+1/\beta)\alpha$-embedding of $Y_{z_i}$ in $H_{z_i}$, it
follows that $\Delta(z_i) = \diam(H_{z_i}) \leq (1+1/\beta)\alpha
\diam(Y_{z_i}) \le (1+1/\beta)\alpha \diam(N_{z_i})$. Let $p$ be a
parent of $U_i$ in $T_M$. Since we have a non-contractive
$\alpha$-embedding of $\U$ in $H_M$ it follows that
$\Delta_{T_M}(p) \ge d_M(A,B)$ for some $A,B \in \U$. Therefore
$\Delta(p) \ge (\beta+1)\gamma \cdot \min \{ d_M(x,y); x \neq y
\in M \}$. Consequently, $\Delta(p)/\Delta(z) \ge
(\beta+1)/[(1+1/\beta)\alpha] \ge k$, by our restriction on
$\beta$. Since $H_M$ and $H_{z_i}$ are $k$-HSTs, it follows that
$T$ also defines a $k$-HST.

It is left to show that $\V$ is $\alpha$-equivalent to $H$. Recall
that for each $z \in M$ there is a non-contractive Lipschitz
bijection $f_z:\U_z\to H_z$ that satisfies for every $A,B \in
\U_z$, $d_{\U_z}(A,B) \le d_{H_z}(f_z(A),f_z(B)) \le \alpha
d_{\U_z}(A,B)$. Define $f:\V\to H$ by $f(V_j^{z_i}) =
f_{z_i}(U_j^{z_i})$. Then, by Claim~\ref{claim:metriccomp} for
every $1\le p<q\le s(z_i)$:
\begin{eqnarray*}
d_\V(V_p^{z_i},V_q^{z_i}) &=& d_{\U_{z_i}}(U_p^{z_i},U_q^{z_i})\\
&\le& d_{H_{z_i}}(f_{z_i}(U_p^{z_i}),f_{z_i}(U_q^{z_i})) =
d_{H}(f(V_p^{z_i}),f(V_q^{z_i}))\\ &\le& (1+1/\beta)\alpha
d_{\U_{z_i}}(U_p^{z_i},U_q^{z_i}) = (1+1/\beta)\alpha
d_\V(V_p^{z_i},V_q^{z_i}) .
\end{eqnarray*}

Additionally, we have a non-contractive Lipschitz bijection
$f_M:\U\to H_M$ that satisfies for every $U_i,U_j \in \U$,
$d_{\U}(U_i,U_j) \le d_{H_M}(f_M(U_i),f_M(U_j)) \le \alpha
d_{\U}(U_i,U_j)$. Hence, by Claim~\ref{claim:metriccomp}, for
every $1\le i<j\le s$ and every $1\le p\le s(z_i)$, $1\le q\le
s(z_j)$:
\begin{eqnarray*}
d_\V(V_p^{z_i},V_q^{z_i}) &\le& (\beta+1)\gamma d_{\U}(U_i,U_j)\\
&\le& (\beta+1)\gamma d_{H_M}(f_M(U_i),f_M(U_j)) =
d_{H}(f(V_p^{z_i}),f(V_q^{z_j}))\\ &\le& \alpha(\beta+1)\gamma
d_{\U}(U_i,U_j) \le  (1+1/\beta)\alpha d_\V(V_p^{z_i},V_q^{z_i}) .
\end{eqnarray*}

The proof of Lemma~\ref{lem:composing} is complete.

\end{proof}

We will also require the following two results from from
\cite{blmn1}:

\begin{lemma}[\cite{blmn1}]\label{lem:compaspect} For any
$\alpha,\beta\ge 1$, if a metric space $M$ is $\alpha$-equivalent
to a $\alpha\beta$-HST, then $M$ is $(1+2/\beta)$-equivalent to a
metric space in $\comp_\beta(\M(\alpha))$.
\end{lemma}

\begin{theorem}[\cite{blmn1}]\label{thm:kHST} There exists a
universal constant $c>0$ such that for every $0<\e\le 1$ and $k\ge
1$ every $n$-point metric space $M$ contains a subset $N\subseteq
M$ which $(2+\e)$-embeds into a $k$-HST and:
$$
|N|\ge n^{\frac{c\e}{\log(2k/\e)}}.
$$

\end{theorem}

We are now in position to present the announced lower bound for
the $QS$ problem for small distortion:

\begin{proposition}\label{prop:QSsmall} There exists a universal constant $C>0$ such that whenever
$M$ is an $n$-point metric space and $0<\e \le 1/2$, there is a
$QS$ space of $M$, $\U$, which is $(1+\e)$-equivalent to a
$1/\e$-HST and:
$$
|\U|\ge n^{\frac{C\e}{[\log(1/\e)]^2}}.
$$
In particular, for every $1\le p\le \infty$:
$$
\QS_p(1+\e,n)\ge n^{\frac{C\e}{[\log(1/\e)]^2}}.
$$
\end{proposition}

\begin{proof} Fix $k\ge 8$  which will be specified later. By
Theorem~\ref{thm:kHST}, $M$ contains a subset $N$ which is
$4$-equivalent to a $k$-HST and $|N|\ge n^{c/\log(2k)}$.
By Lemma~\ref{lem:compaspect}, $N$ is $(1+8/k)$-equivalent to a metric
space in $\comp_{k/4}(\M(4))$. By Corollary~\ref{coro:weightedaspect},
$$
\sigma_{\mathrm{EQ}}\left(M(4),1+\frac{1}{k}\right)\ge
\frac{c'}{k\log k},
$$
for some absolute constant $c'$. By Lemma~\ref{lem:composing},
\begin{eqnarray*}
\sigma_{k/8}\left(\mathrm{comp}_{k/4}(\M(4)),\left(1+\frac{4}{k}\right)\left(1+\frac{1}{k}\right)\right)&\ge&
\sigma_{k/8}\left(\M(4),1+\frac{1}{k}\right)\\&\ge&
\sigma_{\mathrm{EQ}}\left(\M(4),1+\frac{1}{k}\right)\ge
\frac{c'}{k\log k}.
\end{eqnarray*}
Since $N$ is $(1+8/k)$-equivalent to a
 metric
space in $\comp_{k/4}(\M(4))$, it follows that it has a $QS$ space
$\U$ which is $(1+8/k)(1+4/k)(1+1/k)\le1+ 20/k$ equivalent to a
$k/8$-HST, and:
$$
|\U|\ge |N|^{c'/(k\log k)}\ge n^{c''/[k(\log k)^2]},
$$
where $c''$ is an absolute constant. Taking $k=20/\e$ concludes
the proof.
\end{proof}

\section{The QS Problem for the Hypercube}\label{section:cube}

For every integer $d\ge 1$ denote $\Omega_d=\{0,1\}^d$, equipped
with the Hamming ($\ell_1$) metric. Our goal in this section is to
prove Theorem~\ref{thm:cube}, stated in the introduction. As
proved by P. Enflo in \cite{enflocube}, for $1\le p\le 2$,
$c_p(\Omega_d)=d^{1-1/p}$. For $2\le p<\infty$ it was shown in
\cite{gidandi} that there is a constant $a(p)>0$ such that for all
$d$, $c_p(\Omega_d)\ge a(p)\sqrt{d}$. The following lemma
complements these lower bounds:

\begin{lemma}\label{lem:uppercube} For every $1\le p<\infty$ there is an absolute constant
$c=c(p)>0$ such that for every integer $d\ge 1$ and every
$2^{-d}\le\e<1/4$, if $\U$ is a $QS$ space of $\Omega_d$ such that
$|\U|> (1-\e)2^d$ then
$$
c_p(\U)\ge
c\left[\frac{\log(1/\e)}{1+\log\left(\frac{d}{\log(1/\e)}\right)}\right]^{\min\left\{1-\frac{1}{p},\frac12\right\}}.
$$

\end{lemma}

\begin{proof} By adjusting the value of $c$, we may assume that that
$e^{-d/70}< \e<d^{-50}$. In this case, if we set:
$$
r=\frac{1}{16}\left\lfloor
\frac{\log(1/\e)}{\log\left(\frac{d}{\log(1/\e)}\right)}\right\rfloor,
$$
then $3\le r<\frac{d}{16}$. The ball of radius $2r$ in $\Omega_d$
contains at most  $3r\binom{d}{2r}\le
3r\left(\frac{ed}{2r}\right)^{2r}\le e^{4r\log(d/r)}$ points.
Therefore, the cube $\Omega_d$ contains at least $2^{d}\cdot
e^{-4r\log(d/r)}$ disjoint balls of radius $r$. Writing
$x=d/\log(1/\e)$ we have $16x\log x\le x^2$, so that $4r\log
(d/r)\le \frac{1}{4}\log(1/\e)\cdot\frac{\log[16x\log x]}{\log x}
\le \log[1/(2\e)]$. This reasoning shows that $\Omega_d$ contains
at least $2\e 2^d$ disjoint balls of radius $r$.

Let $\U=\{U_1,\ldots,U_k\}$ be a $QS$ space of $\Omega_d$ with $k>
(1-\e)2^d$. As in the proof of Lemma~\ref{lem:upperexpander}, $\U$
must contain more than $(1-2\e)2^d$ singletons. Since $\Omega_d$
contains at least $2\e2^d$ disjoint balls of radius $r$, it
follows that $\U$ must contain the elements of some ball $B$ of
radius $r$ as singletons. Let $x$ be the center of $B$. Write
$m=\left\lfloor \frac{r}{3}\right\rfloor$ and consider the
sub-cube $C=\{0,1\}^m\times\{x_{m+1}\}\times\cdots\times\{x_d\}$.
Observe that $C\subseteq B$, and the diameter of $C$ is at most
$2r/3$. Moreover, since $\U$ contains the elements of $B$ as
singletons, the distance in $\Omega_d$ between an element of $C$
and a non-singleton element of $\U$ is at least $2r/3$. This shows
that when calculating the geodesic distance in $\U$ between two
points in $C$, it is enough to restrict ourselves to paths which
pass only through singletons. It follows that the metric induced
by $\U$ on $C$ coincides with the Hamming metric. By the results
of \cite{enflocube} and \cite{gidandi}, it follows that
$c_p(\U)\ge c_p(C)\ge a(p) k^{\min\{1-1/p,1/2\}}$, for some
constant $a(p)$ depending only on $p$. This completes the proof.

\end{proof}

We now turn our attention to the construction of large $QS$ spaces
of the hypercube which well embed into $\ell_p$. Our proof yields
several embedding results which may be useful in other
circumstances. The case $p=2$ is simpler, so deal with it first.

Given a metric space $M$ and $D>0$, we denote by $M^{\le D}$ the
metric space $(M,d_{M^{\le D}})$, where $d_{M^{\le
D}}(x,y)=\min\{d_M(x,y),D\}$.

\begin{lemma}\label{lem:truncation} For every $D>0$,
$ c_2\big(\ell_2^{\le D}\big)\le \sqrt{\frac{e}{e-1}} $. In fact,
$\ell_2^{\le D}$ $\sqrt{\frac{e}{e-1}} $-embeds into the
$\ell_2$-sphere of radius $D$.
\end{lemma}
\begin{proof}
Let $\{g_i\}_{i=1}^\infty$ i.i.d. standard Gaussian random
variables. Assume that they are defined on some probability space
$\Omega$. Consider the Hilbert space $H=L_2(\Omega)$ where we
think of $L_2(\Omega)$ as all the complex valued square integrable
functions on $\Omega$. Define $F:\ell_2\to H$ by:
$$
F(x_1,x_2,\ldots)=D\exp\left(\frac{i}{D}\sum_{j=1}^\infty
x_jg_j\right).
$$
Clearly $\|F(x)\|_2=D$ for every $x\in \ell_2$. Observe that for
every $x,y\in \ell_2$,
\begin{eqnarray*}
|F(x)-F(y)|^2&=&\\&=&
D^2\left|\exp\left(\frac{i}{D}\sum_{j=1}^\infty
x_jg_j\right)-\exp\left(\frac{i}{D}\sum_{j=1}^\infty
y_jg_j\right)\right|^2\\
&=&D^2\left|\exp\left(\frac{i}{D}\sum_{j=1}^\infty
y_jg_j\right)\left[\exp\left(\frac{i}{D}\sum_{j=1}^\infty
(x_j-y_j)g_j\right)-1\right]\right|^2\\
&=&D^2\left|\exp\left(\frac{i}{D}\sum_{j=1}^\infty
(x_j-y_j)g_j\right)-1\right|^2\\
&=&2D^2\left[1-\cos\left(\frac{1}{D}\sum_{j=1}^\infty
(x_j-y_j)g_j\right)\right].
\end{eqnarray*}

Now, $\sum_{i=1}^\infty (x_j-y_j)g_j$ has the same distribution as
$g_1\sqrt{\sum_{j=1}^\infty (x_j-y_j)^2}$. Hence:
\begin{eqnarray*}
\E|F(x)-F(y)|^2=2D^2\left[1-\E\cos\left(\frac{g_1}{D}\|x-y\|_2\right)\right].
\end{eqnarray*}
Observe that by symmetry,
$\E\sin\left(\frac{g_1}{D}\|x-y\|_2\right)=0$, so that:
\begin{eqnarray*}
\E\cos\left(\frac{g_1}{D}\|x-y\|_2\right)=
 \E\exp\left(i\frac{g_1}{D}\|x-y\|_2\right)
= \exp\left(-\frac{\|x-y\|_2^2}{2D^2}\right),
\end{eqnarray*}
where we use the fact that $\E e^{iag_1}=e^{-a^2/2}$.

Putting it all together, we have shown that:
$$
\|F(x)-F(y)\|_2=\sqrt{2}D\sqrt{1-e^{-\frac{\|x-y\|_2^2}{2D^2}}}.
$$

Using the elementary inequality:
$$
\frac{e-1}{e}\min\{1,a\}\le 1-e^{-a}\le \min\{1,a\} \quad a>0,
$$
we deduce that:
$$
\sqrt{\frac{e-1}{e}}\min\{D,\|x-y\|_2\}\le\|F(x)-F(y)\|_2\le
\min\{D,{\|x-y\|_2}\}.
$$
\end{proof}

It should be remarked here that in the proof of
Lemma~\ref{lem:truncation} we have only used the fact that the
metric $\sqrt{1-e^{-\|x-y\|^2_2}}$, defined on $\ell_2$, embeds
isometrically in Hilbert space. This follows from a classical
result of Schoenberg~\cite{schoenberg}, but we have chosen to
present the above direct proof since, apart from giving a concrete
embedding, this argument generalizes to $\ell_p$, $1\le p\le 2$
(Lemma~\ref{lem:ptruncation} below).

\begin{remark} \label{rem:1}
 Lemma~\ref{lem:truncation} cannot
be replaced by an isometric result. In fact, for every $D>0$,
$$
c_2(\ell_2^{\le D})\ge \frac{2\sqrt{5-\sqrt{7}}}{3}>1.02.
$$
To see this let $T:\R^2\to \ell_2$ be such that for every $x,y\in
\R^2$, $\min\{\|x-y\|_2,D\}\le \|T(x)-T(y)\|_2\le
A\min\{\|x-y\|_2,D\}$. It is straightforward to verify that when
viewed as a subset of $\ell_2^{\le D}$, the points
$\{(0,0),(D,0),(D/2,D),(D/2,0)\}$ cannot be isometrically embedded
in Hilbert space. To lower-bound the distortion, define
$a=T(0,0)$, $b=T(D,0)$, $c=T(D/2,D)$, $d=T(D/2,0)$. By the
parallelogram identity:
$$
\frac{D^2A^2}{2}\ge
\|a-d\|_2^2+\|b-d\|_2^2=\frac{\|a+b-2d\|_2^2+\|a-b\|_2^2}{2}\ge
2\left\|\frac{a+b}{2}-d\right\|_2^2+\frac{D^2}{2}.
$$
Hence:
$$
\left\|\frac{a+b}{2}-d\right\|_2\le \frac{D\sqrt{A^2-1}}{2}.
$$
Similarly:
$$
2D^2A^2\ge \|a-c\|_2^2+\|b-c\|_2^2=
\frac{\|a+b-2c\|_2^2+\|a-b\|_2^2}{2}\ge
2\left\|\frac{a+b}{2}-c\right\|_2^2+\frac{D^2}{2},
$$
or
$$
\left\|\frac{a+b}{2}-c\right\|_2\le \frac{D\sqrt{4A^2-1}}{2}.
$$
But:
$$
D\le \|c-d\|_2\le
\left\|\frac{a+b}{2}-d\right\|_2+\left\|\frac{a+b}{2}-c\right\|_2\le
\frac{D\sqrt{A^2-1}}{2}+\frac{D\sqrt{4A^2-1}}{2},
$$
which simplifies to give the required result.
\end{remark}

\begin{remark} \label{rem:2}
Let $\omega:[0,\infty)\to [0,\infty)$ be a
concave non-decreasing function such that $\omega(0)=0$ and
$\omega(t)>0$ for $t>0$. It is straightforward to verify that if
we define for $x,y\in \ell_2$, $d_\omega(x,y)=\omega(\|x-y\|_2)$,
then $d_\omega$ is a metric. Lemma~\ref{lem:truncation} dealt with
the case $\omega(t)=\min\{t,D\}$, but we claim that in fact there
is a constant $C>0$ such that for every such $\omega$,
$c_2(\ell_2,d_\omega)\le C$. To see this observe that
$\phi(t)=\omega\left(\sqrt{t}\right)^2$ is still concave and
non-decreasing (assume by approximation that $\omega$ is
differentiable and observe that
$\phi'(t^2)=\omega'(t)\omega(t)/t$. By our assumptions, both
$\omega'(t)$ and $\omega(t)/t$ are non-negative and
non-increasing, so that the required result follows). Now, it is
well known (see for example Proposition 3.2.6. in~\cite{brudnyi})
that we may therefore write $\phi(t)\approx \sum_{i=1}^\infty
\min\{\lambda_i,\mu_i t\}$ for some $\lambda_i,\mu_i>0$, where the
symbol $\approx$ means that the two functions are equivalent up to
absolute multiplicative constants. By Lemma~\ref{lem:truncation},
for every $i$ there is a function $F_i:\ell_2\to \ell_2$ such that
for every $x,y\in \ell_2$ $\|F_i(x)-F_i(y)\|_2\approx
\min\{\sqrt{\lambda_i},\sqrt{\mu_i}\|x-y\|_2\}$. Define
$F:\ell_2\to \ell_2(L_2)$ by setting the $i$'th coordinate of $F$
to be $F_i$. Then for every $x,y\in \ell_2$:
$$
\|F(x)-F(y)\|_2^2=\sum_{i=1}^\infty \|F_i(x)-F_i(y)\|_2^2\approx
\sum_{i=1}^\infty \min\{\lambda_i,\mu_i\|x-y\|_2^2\}\approx
\phi(\|x-y\|_2^2)=\omega(\|x-y\|_2)^2.
$$
\end{remark}

\begin{lemma}\label{lem:sqrt} Let $X$ be a metric space such that
$\min_{x\neq y}d_X(x,y)\ge 1$ and the metric space
$(X,\sqrt{d_X})$ is isometric to a subset of $\ell_2$. Then for
every $D\ge 1$, $c_2\big(X^{\le D}\big)\le \sqrt{\frac{eD}{e-1}}$.
Moreover, exists a $1$-Lipschitz embedding $f:X^{\le D}\to \ell_2$
such that $\mathrm{dist}(f)\le \sqrt{\frac{eD}{e-1}}$ and for
every $x\in X$, $\|f(x)\|_2= \sqrt{D}$.
\end{lemma}

\begin{proof} All we have to do is to observe that for every
$x,y\in X$,
$$
\min\{\sqrt{D},\sqrt{d_X(x,y)}\}\le \min\{D,d_X(x,y)\}\le
\sqrt{D}\cdot \min\{\sqrt{D},\sqrt{d_X(x,y)}\},
$$
and then apply Lemma~\ref{lem:truncation}.

\end{proof}

\begin{corollary}\label{coro:truncatedcube} For every integer
$d\ge 1$, $c_2\big(\Omega_d^{\le D}\big)\le
\sqrt{\frac{eD}{e-1}}$, where the embedding is $1$-Lipschitz and
takes values in the $\ell_2$-sphere of radius $\sqrt{D}$.
\end{corollary}

\begin{proof} This follows from Lemma~\ref{lem:sqrt} and the
classical fact \cite{ww} that $\ell_1$ equipped with the metric
$\sqrt{\|x-y\|_1}$ is isometric to a subset of $\ell_2$.

\end{proof}

\begin{lemma}\label{lem:largeQScube}
There is a universal constant $C>0$ such that for every integer
$d\ge 1$ and every $2^{-d}\le \e<1/4$ there exists a $QS$ space of
$\Omega_d$, $\U$, such that $|\U|\ge (1-\e)2^d$ and:
$$
c_2(\U)\le
C\sqrt{\frac{\log(1/\e)}{1+\log\left(\frac{d}{\log(1/\e)}\right)}}.
$$
\end{lemma}

\begin{proof}
By adjusting the constant $C$, we may assume that $\e\ge
e^{-d/400}$. Define $r$ to be the smallest even integer greater
than:
$$
2\left\lceil\frac{\log (1/\e)}{\log
\left(\frac{d}{\log(1/\e)}\right)}\right\rceil.
$$
We first construct a subset $A\subseteq \Omega_d$ via the
following iterative procedure: Pick any $x_1\in \Omega_d$. Having
chosen $x_1,\ldots,x_{k-1}$, as long as $\Omega_d\setminus
\cup_{j=1}^{k-1}B_{\Omega_d}(x_j,2r)\neq \emptyset$, pick any
$x_i\in \Omega_d\setminus \cup_{j=1}^{k-1}B_{\Omega_d}(x_j,2r)$.
When this procedure terminates we set $A=\{x_1,x_2,\ldots\}$.

Define $S\subseteq \Omega_d$ by:
$$
S=\Omega_d\setminus \left(\bigcup_{x\in A}
B_{\Omega_d}(x,r/2)\setminus\{x\}\right).
$$
The $QS$ space of $\Omega_d$ which we consider is $\U=S/A$.

We first bound the cardinality of $\U$ from below. Observe that by
the construction, the balls $\{B_{\Omega_d}(x,r)\}_{x\in A}$ are
disjoint, so that $|A|\binom{d}{r}\le 2^d$. Hence:
$$
|\U|\ge 2^d-|A|\cdot
r\binom{d}{r/2}+1>\left(1-\frac{r\binom{d}{r/2}}{\binom{d}{r}}\right)2^d\ge
\left(1-\frac{r\left(\frac{ed}{r/2}\right)^{r/2}}{\left(\frac{d}{r}\right)^r}\right)2^d=
\left(1-re^{-\frac{r}{2}\log\left(\frac{d}{2er}\right)}\right)2^d.
$$
By our choice of $r$, and the restriction $\e\ge  e^{-d/400}$, it
is straightforward to verify that
$re^{-\frac{r}{2}\log\left(\frac{d}{2er}\right)}\le \e$. We have
shown that $|\U|\ge (1-\e)2^d$.

By our construction, for every $x\in S\setminus A$, $r/2\le
d_{\Omega_d}(x,A)\le 2r$. This implies that for every $x,y\in
\U\setminus \{A\}$,
$$\min\{d_{\Omega_d}(x,y),r\} \le d_\U(x,y)\le \min\{d_{\Omega_d}(x,y),4r\}.
$$
By Corollary~\ref{coro:truncatedcube} there is an embedding
$f:\U\setminus \{A\}\to \ell_2$ such that for every $x\in
\U\setminus \{A\}$, $\|f(x)\|_2=\sqrt{r}$ and for every $x,y\in
\U\setminus \{A\}$:
$$
\sqrt{\frac{e-1}{16er}}\cdot d_\U(x,y)\le \|f(x)-f(y)\|_2\le
d_\U(x,y).
$$
Since for every $x\in \U\setminus \{A\}$, $r/2\le d_\U(x,A)\le
2r$, we may extend $f$ to $\U$ by setting $f(A)=0$. As $f$ takes
values in the $\ell_2$-sphere of radius $\sqrt{r}$,
$\mathrm{dist}(f)=O(\sqrt{r})$, as required.
\end{proof}

Since $\ell_2$ embeds isometrically into $L_p$, $p\ge 1$,
Lemma~\ref{lem:largeQScube} implies that Lemma~\ref{lem:uppercube}
is optimal (up to the dependence of the constant on $p$) for $p\ge
2$. The case $1\le p\le 2$ seems to be more delicate, but we can
still match the bound in Lemma~\ref{lem:uppercube} up to
logarithmic factors.

Recall that for $1\le p\le 2$ there exists a symmetric $p$-stable
random variable $g$. This means that there exists a constant
$c=c(p)>0$ such that for every $t\in \R$, $\E
e^{itg}=e^{-c|t|^p}$. In what follows we fix $1\le p<2$ and ignore
the dependence of all the constants on $p$. Moreover, given two
quantities $A,B$ the notation $A\approx_p B$ means that there are
constants $C_1,C_2$, which may depend only on $p$, such that $C_1
A\le B\le C_2 A$. Denote the density of $g$ by $\varphi$. It is
well known (see \cite{zolotarev}) that $ \varphi(t)\approx_p
\frac{1}{1+t^{p+1}}$.

\begin{lemma}\label{lem:cos} Fix $1\le p<2$ and let $g$ be a
symmetric $p$-stable random variable. Then for every $a>0$,
$$
\E\left[1-\cos(ag)\right]^{p/2}\approx_p
\min\left\{a^p\log\left(\frac{1}{a}+1\right),1\right\}.
$$
\end{lemma}

\begin{proof} Since for $0\le x\le 1$, $1-\cos x\approx x^2$, we have
that:
\begin{eqnarray*}
\E\left[1-\cos(ag)\right]^{p/2}&=&2\int_0^\infty
\left[1-\cos(au)\right]^{p/2}\varphi(u)du\\
&\approx_p& \int_0^{1/a}\frac{a^p
u^p}{1+u^{p+1}}du+\int_{1/a}^\infty \frac{1}{1+u^{p+1}}du\approx_p
\min\left\{a^p\log\left(\frac{1}{a}+1\right),1\right\}.
\end{eqnarray*}
\end{proof}

The following lemma is analogous to Lemma~\ref{lem:truncation}:

\begin{lemma}\label{lem:ptruncation} For every $1\le p<2$ and every $D>0$ there
exists a mapping $F:\ell_p\to L_p$ such that for every $x\in
\ell_p$, $\|F(x)\|_p=D$ and for every $x,y\in \ell_p$,
$$
\|F(x)-F(y)\|_p\approx_p
\min\left\{\|x-y\|_p\left[\log\left(\frac{D}{\|x-y\|_p}+1\right)\right]^{1/p},D\right\}.
$$
\end{lemma}

\begin{proof}
Let $\{g_i\}_{i=1}^\infty$ i.i.d. symmetric $p$-stable random
variables. Assume that they are defined on some probability space
$\Omega$. Consider the space $L_p(\Omega)$, where we think of
$L_p(\Omega)$ as all the complex valued $p$-integrable functions
on $\Omega$. Define $F:\ell_p\to H$ by:
$$
F(x_1,x_2,\ldots)=D\exp\left(\frac{i}{D}\sum_{j=1}^\infty
x_jg_j\right).
$$
Clearly $\|F(x)\|_p=D$ for every $x\in \ell_p$. As we have seen in
the proof of Lemma~\ref{lem:truncation}, for $x,y\in \ell_p$:
\begin{eqnarray*}
|F(x)-F(y)|^p=2^pD^p\left[1-\cos\left(\frac{1}{D}\sum_{j=1}^\infty
(x_j-y_j)g_j\right)\right]^{p/2}.
\end{eqnarray*}

Now, $\sum_{i=1}^\infty (x_j-y_j)g_j$ has the same distribution as
$g_1\|x-y\|_p$. Hence by Lemma~\ref{lem:cos}:
\begin{eqnarray*}
\E|F(x)-F(y)|^p&=&2^pD^p\E\left[1-\cos\left(\frac{g_1}{D}\|x-y\|_p\right)\right]^{p/2}\\&\approx_p&
D^p\min\left\{\frac{\|x-y\|_p^p}{D^p}\log\left(\frac{D}{\|x-y\|_p}+1\right),1\right\}.
\end{eqnarray*}
\end{proof}

\begin{remark} \label{rem:3}
The above argument also shows that for every $1\le q<p\le 2$ there
is a constant $C=C(p,q)$ such that for every $D>0$, $\ell_p^{\le
D}$ is $C$-equivalent to a subset of $L_q$ (since in this case
there is no logarithmic term in Lemma~\ref{lem:cos}). For every
$1\le q<p\le 2$, the metric space $(L_q,\|x-y\|_q^{q/p})$ is
isometric to a subset of $L_p$. When $p\le 2$ this follows from
general results of Bretagnolle, Dacunha-Castelle and
Krivine~\cite{krivine} (see also the book~\cite{ww}). It is of
interest, however, to give a concrete formula for this embedding,
which works for every $1\le q<  p<\infty$. To this end observe
that by a change of variable it follows that for every
$0<\alpha<2\beta$ there exists a constant $c_{\alpha,\beta}>0$
such that for every $x\in \R$, $
|x|^\alpha=c_{\alpha,\beta}\int_{-\infty}^\infty \frac{(1-\cos
tx)^\beta}{|t|^{\alpha+1}} dt $. Define $T: L_q(\R)\to
L_p(\R\times \R)$ by
$T(f)(s,t)=\frac{1-e^{itf(s)}}{|t|^{(q+1)/p}}$. For every $f,g\in
L_p(\R)$ we have:
\begin{eqnarray*}
\|T(f)-T(g)\|_p^p&=&\int_{-\infty}^\infty \int_{-\infty}^\infty
\frac{|1-e^{it[f(s)-g(s)]}|^p}{|t|^{q+1}}dt
ds\\&=&2^{p/2}\int_{-\infty}^\infty \int_{-\infty}^\infty
\frac{\{1-\cos[t(f(s)-g(s))]\}^{p/2}}{|t|^{q+1}}dtds=
2^{p/2}c_{q+1,p/2}\|f-g\|_q^q,
\end{eqnarray*}
so that $T$ is the required isometry.

A corollary of these observations is that for every $\e>0$ there
is a constant $C(\e)>0$ such that for every $D>0$ the metric
$\min\{\|x-y\|_p^{1-\e},D\}$ on $\ell_p$, $1\le p\le 2$, is
$C(\e)$-equivalent to a subset of $L_p$. We do not know whether
the exponent $1-\e$ can be removed in this statement.
\end{remark}

\begin{remark}
The same reasoning as in Remark~\ref{rem:2} shows that for every
$\omega:[0,\infty)\to [0,\infty)$ which is concave,
non-decreasing, $\omega(0)=0$ and $\omega(t)>0$ for $t>0$, the
metric $\omega(\|x-y\|_p)$ on $\ell_p$ is $c(p,q)$-equivalent to a
subset of $L_q$ for every $1\le q<p\le 2$. The only difference in
the proof is that one should apply the same argument to show that
$\phi(t)=\omega(t^{1/q})^q$ shares the same properties as
$\omega$. Similarly, the metric $\omega(\|x-y\|_p^{1-\e})$ is
$C(\e)$-equivalent to a subset of $L_p$.
\end{remark}

\begin{remark}
Remark~\ref{rem:3} is false for $p>2$. In
fact, for every $0<\gamma\le 1$, and $D>0$, the metric
$\min\{\|x-y\|_p^\gamma,D\}$ is not Lipschitz equivalent to a
subset of $L_q$ for any $1\le q<\infty$. To see this observe that
if we assume the contrary then this metric would be Lipschitz
equivalent to a bounded subset of $L_q$. An application of Mazur's
map (see~\cite{benlin}) shows that this implies that $L_p$ is
uniformly homeomorphic to a subset of $L_2$. Since $L_p$, $p>2$
has type 2, a theorem of Aharoni, Maurey and
Mityagin~\cite{aharoni} implies that $L_p$ would be linearly
isomorphic to a subspace of $L_1$. This is a contradiction since
$L_1$ has cotype 2 while $L_p$, $p>2$ has cotype $p$. Actually, by
the results presented in Chapter 9 of~\cite{benlin}, this argument
implies that the above metric is not Lipschitz equivalent to a
subset of any separable Banach lattice with finite cotype.
\end{remark}

\medskip

\begin{corollary}\label{coro:uptolog} Fix $1\le p<2$. Let $X$ be a
finite subset of $L_1$ such that for every $x,y\in X$,
$\|x-y\|_1\ge 1$. Then for every $D\ge 2$ there is an embedding
$\psi:X\to L_p$ such that for every $x\in X$, $
\|\psi(x)\|_p=D^{1/p}$ and for every $x,y\in X$,
$$
\frac{C_1(p)}{D^{1-1/p}}\min\{\|x-y\|_1,
D\}\le\|\psi(x)-\psi(y)\|_p\le C_2(p)(\log
D)^{1/p}\min\{\|x-y\|_1, D\},
$$
where $C_1(p),C_2(p)$ are constants which depend only on $p$.
\end{corollary}

\begin{proof} We begin by noting that as in Remark~\ref{rem:3}, there is a mapping $G:L_1\to
L_p$ such that for every $x,y\in L_1$,
$\|G(x)-G(y)\|_p=\|x-y\|_1^{1/p}$. Since $X$ is finite, there is
an isometric embedding $T:G(X)\to \ell_p$ (see~\cite{dezalau}).
Let $F$ be as in Lemma~\ref{lem:ptruncation}, with $D$ replaced by
$D^{1/p}$, and define $\psi=F\circ T\circ G$. Now,
$\|\psi(x)\|_p=D^{1/p}$ for every $x\in X$ and:
\begin{eqnarray*}
\|\psi(x)-\psi(y)\|_p&\approx_p&
\min\left\{\|x-y\|_1^{1/p}\cdot\left[\log\left(\frac{D}{\|x-y\|_1}+1\right)\right]^{1/p},
D^{1/p}\right\}\\
&\le& C(p)(\log D)^{1/p}\min\{\|x-y\|_1^{1/p},D^{1/p}\}\\&\le&
C(p)(\log D)^{1/p}\min\{\|x-y\|_1,D\},
\end{eqnarray*}
where we have used the fact that $\|x-y\|_1\ge 1$. Similarly, we
have the inequality:
$$
\|\psi(x)-\psi(y)\|_p\ge
\frac{C'(p)}{D^{1-1/p}}\min\{\|x-y\|_1,D\}.
$$
The proof is complete.
\end{proof}

A proof identical to the proof of Lemma~\ref{lem:largeQScube} now
gives a bound which nearly matches the bound in
Lemma~\ref{lem:uppercube}:

\begin{lemma}\label{lem:plargeQScube} For every $1\le p<2$ there is a constant $C(p)>0$
such that for every integer $d\ge 1$ and every $2^{-d}\le \e<1/4$
there exists a $QS$ space of $\Omega_d$, $\U$, such that $|\U|\ge
(1-\e)2^d$ and:
$$
c_p(\U)\le
C(p)\left[\frac{\log(1/\e)}{1+\log\left(\frac{d}{\log(1/\e)}\right)}\right]^{1-1/p}\cdot\left[\log\left(
\frac{\log(1/\e)}{1+\log\left(\frac{d}{\log(1/\e)}\right)}\right)\right]^{1/p}.
$$

\end{lemma}

\section{Lipschitz Quotients}\label{section:lipq}

In this section we prove Theorem~\ref{thm:lipqs}. We shall use the
notation introduced in the introduction.

Recall that for a metric space $X$ and two subsets $U,V\subset X$,
the Hausdorff distance between $U$ and $V$ is defined as:
$$
\HH_X(U,V)=\sup\{\max\{d_X(u,V),d_X(v,U)\};\ u\in U,\ v\in V\}.
$$

The following straightforward lemma is the way we will use the
Lipschitz and co-Lipschitz conditions:

\begin{lemma}\label{lem:lipq} Let $X$, $Y$ be metric spaces and $A>0$. For
every surjection $f:X\to Y$ the following assertions hold:
\begin{enumerate}
\item $\Lip(f)\le A$ if and only if for every $y,z\in Y$,
$d_Y(y,z)\le Ad_X(f^{-1}(y),f^{-1}(z))$. \item $\coLip(f)\le A$ if
and only if for every $y,z\in Y$, $\HH_X(f^{-1}(y),f^{-1}(z))\le A
d_Y(y,z)$.
\end{enumerate}

\end{lemma}

\begin{remark}\label{rem:singletons} A simple corollary of Lemma~\ref{lem:lipq}, which
will be useful later, is that if $f:X\to Y$ is a Lipschitz
quotient and we set $U=f^{-1}(\{y\in Y;\ |f^{-1}(y)|=1\})$ then
$f|_U$ is a Lipschitz equivalence between $U$ and $f(U)$ and
$\mathrm{dist}(f|_U)\le \Lip(f)\cdot\coLip(f)$.

\end{remark}

\medskip

In the following lemma we prove two recursive inequalities which
will be used to give upper bounds for $\QS^{\Lip}$. In this lemma
we use the notation $R_\M(\cdot,\cdot)$ which was introduced in
the introduction.

\begin{lemma}\label{lem:lipcomp} Let $\M$ be a class of metric spaces. Then for
every two integers $k,m\ge 1$ and every $\alpha\ge 1$,
\begin{enumerate}
\item
$\QS_\M^{\Lip}(\alpha,km)\le\QS_\M^{\Lip}(\alpha,k)\cdot\QS_\M^{\Lip}(\alpha,m)$.
\item $\QS_\M^{\Lip}(\alpha,km)\le k+R_\M(\alpha,k)\cdot
\QS_\M^{\Lip}(\alpha,m)$.
\end{enumerate}
\end{lemma}

\begin{proof} We will start by proving the first assertion.
Denote $a=\QS_\M^{\Lip}(\alpha,k)$, $b=\QS_\M^{\Lip}(\alpha,m)$.
Let $X$ be a $k$-point metric space such that the largest
$\alpha$-Lipschitz quotient of a subspace of $X$ which is in a
member of $\M$ has $a$ points. Similarly, let $Y$ be an $m$-point
metric space such that the largest $\alpha$-Lipschitz quotient of
a subspace of $Y$ which is in a member of $\M$ has $b$ points. We
think of $X$ as a metric $d_X$ on $[k]=\{1,\ldots,k\}$. Fix any:
$$
\mu>\alpha\Phi(Y)=\frac{\alpha\diam(Y)}{\min_{y\neq
z}d_Y(y,z)}\qquad \mathrm{and}\qquad \theta\ge
\alpha\mu^k\frac{\diam(Y)}{\min_{1\le i<j\le k}d_X(i,j)}.
$$

Set $Z=Y\times [k]$ and define:
\[ d_Z((y,i),(z,j))= \begin{cases} \mu^i d_Y(y,z) & i=j \\
  \theta d_X(i,j) & i\neq j .\end{cases}  \]

This definition is clearly a particular case of metric
composition, and the choice of parameters ensures that $d_Z$ is
indeed a metric.

Assume that there is $S\subseteq Z$, $M\in \M$ and $N\subseteq M$
such that there is an $\alpha$-Lipschitz quotient $f:S\to N$. Our
goal is to show that $|N|\le ab$.

Observe that by the definition of the metric on $Z$ we have that
for every $i\in [k]$ and $p,q\in N$, $p\neq q$, if $f^{-1}(p)\cap
(Y\times\{i\}),f^{-1}(q)\cap (Y\times\{i\})\neq\emptyset$ then
$d_Z(f^{-1}(p),f^{-1}(q))\le \mu^i\diam(Y)\le \mu^k\diam(Y)$. On
the other hand, if in addition for some $j\in [k]$, $j\neq i$,
$f^{-1}(p)\cap (Y\times\{j\})\neq \emptyset$ but $f^{-1}(q)\cap
(Y\times\{j\})=\emptyset$ then $\HH_Z(f^{-1}(p),f^{-1}(q))\ge
\theta \min_{1\le i<j\le k}d_X(i,j)>\alpha\mu^k\diam(Y) $. This is
a contradiction since Lemma~\ref{lem:lipq} implies in particular
that
\begin{eqnarray}\label{eq:ratio}
\frac{\HH_Z(f^{-1}(p),f^{-1}(q))}{d_Z(f^{-1}(p),f^{-1}(q))}\le
\alpha.
\end{eqnarray}
Hence $f^{-1}(q)\cap (Y\times\{j\})\neq\emptyset$. Without loss of
generality assume that $j>i$. Then:
$$
\HH_Z(f^{-1}(p),f^{-1}(q))\ge \mu^j\min_{y\neq z}
d_Y(y,z)>\mu^{j-1}\alpha\diam(Y)\ge \mu^i\alpha\diam(Y),
$$
and we arrive once more to a contradiction with (\ref{eq:ratio}).

Summarizing, we have shown that for every $i\in [k]$ and $p\in N$,
if $f^{-1}(p)\cap (Y\times\{i\})\neq \emptyset$ then either
$f^{-1}(p)\subseteq Y\times \{i\}$ or $f^{-1}(p)\supseteq
f^{-1}(N)\cap (Y\times\{i\})$. In particular, if we write for
$p,q\in N$, $p\sim q$ if there is $i\in [k]$ such that
$f^{-1}(p)\cap (Y\times\{i\})\neq \emptyset$ and $ f^{-1}(q)\cap
(Y\times\{i\})\neq \emptyset$. Then $\sim$ is an equivalence
relation. Let $C_1,\ldots,C_s$ be the equivalence classes of
$\sim$. Take any $p_j\in C_j$ and let $A_j\subset [k]$ be the set
of indices $i\in [k]$ such that there exists $y\in Y$ for which
$(y,i)\in f^{-1}(p_j)$. By the definition of $\sim$,
$A_1,\ldots,A_s$ are disjoint. Let $A=\cup_{j=1}^s A_j$ and define
$g:A\to \{p_1,\ldots,p_s\}$ by: if $i\in A_j$ then $g(i)=p_j$. By
the definition of $d_Z$, if $j\neq \ell$ and $h\in A_j$, $i\in
A_\ell$ then for every $y,z\in Y$,
$d_X(h,i)=d_Z((y,h),(z,i))/\theta$. Hence
$d_A(g^{-1}(p_j),g^{-1}(p_\ell))=d_Z(f^{-1}(p_j),f^{-1}(p_\ell))/\theta$,
$\HH_A(g^{-1}(p_j),g^{-1}(p_\ell))=\HH_Z(f^{-1}(p_j),f^{-1}(p_\ell))/\theta$.
By Lemma~\ref{lem:lipq}, $g$ is an $\alpha$-Lipschitz quotient
from the subspace $A\subset [k]$ onto $\{p_1,\ldots, p_s\}$. It
follows that $s\le a$.

We will conclude once we show that for every $j$, $|C_j|\le b$. If
$|C_j|=1$ then there is nothing to prove. Otherwise there is $i\in
[k]$ such that for every $p\in C_j$, $f^{-1}(p)\subset Y\times
\{i\}$. Lemma~\ref{lem:lipq} implies that $f|_{f^{-1}(C_j)}$ is an
$\alpha$-Lipschitz quotient of a subspace of $Y\times\{i\}$, and
since the metric on $Y\times\{i\}$ is a dilation of $d_Y$, it
follows from the definition of $b$ that $|C_j|\le b$.

\medskip

To prove the second assertion in Lemma~\ref{lem:lipcomp} we repeat
the same construction, but now with $Y$ as before, and $X$ a
$k$-point metric space whose largest subspace which
$\alpha$-embeds into a member of $\M$ has $c=R_\M(\alpha,k)$
points. The rest of the notation will be as above.

Consider the equivalence classes $C_1,\ldots,C_s\subseteq N$, and
enumerate them in such a way that $|C_1|=\ldots=|C_t|=1$ and
$|C_{t+1}|,\ldots,|C_s|\ge 2$. As we have seen above, for $1\le
j\le t$, since $C_j=\{p_j\}$, there is a subset $I_j\subset [k]$
such that $f^{-1}(p_j)=f^{-1}(N)\cap(Y\times I_j)$. Since
$I_1,\ldots,I_t$ are disjoint, $t\le k$. Now, by the construction,
for $t<j\le s$, $|g^{-1}(p_j)|=1$, so that by
Remark~\ref{rem:singletons} we get that
$\{g^{-1}(p_{t+1}),\ldots,g^{-1}(p_s)\}\subseteq [k]$
$\alpha$-embeds into $N$. By the definition of $c$, it follows
that $s-t\le c$. Finally, we have also shown that for every $j$
$|C_j|\le b$, so that $|N|\le t+(s-t)b\le k+cb$, as required.
\end{proof}

\begin{corollary}\label{coro:lipqupper} Let $\M$ be a class of
finite metric spaces and $\alpha\ge 1$. Assume that there is a
finite metric space $M$ such that $c_\M(M)>\alpha$. Then there is
$0\le \delta<1$ such that for infinitely many $n$'s,
$\QS_M^{\Lip}(\alpha,n)\le n^\delta$.
\end{corollary}

\begin{proof} Our assumption implies that there is $0\le \delta<1$ such that
$\QS_\M^{\Lip}(\alpha,|M|)\le |M|-1\le |M|^\delta$. An iteration
of Lemma~\ref{lem:lipcomp} now implies that for every $i\ge 1$,
$\QS_\M^{\Lip}(\alpha,|M|^i)\le |M|^{\delta i}$.
\end{proof}

\begin{remark}\label{rem:tightlarge} For every $1\le p<\infty$ and
$\alpha>2$ there is an integer $n_0=n_0(p,\alpha)$ and constants
$c=c(p,\alpha)$, $C=C(p,\alpha)$ such that $0<c\le C<1$ and for
every $n\ge n_0$, $n^c\le \QS_p^{\Lip}(\alpha,n)\le n^C$. This
follows from Corollary~\ref{coro:lipqupper} and the trivial
inequality $\QS_p^{\Lip}(\alpha,n)\ge R_p(\alpha,n)$, together
with the results of \cite{blmn1}.

\end{remark}

\begin{corollary}\label{coro:lipquppersmall} For every $1\le
\alpha<2$ and $1\le p\le 2$ there is an integer
$n_0=n_0(p,\alpha)$ and a constant $C=C(p,\alpha)$ such that for
every $n\ge n_0$:
$$
\QS_p^{\Lip}(\alpha,n)\le e^{C\sqrt{(\log n)(\log \log n)}}.
$$
For $p>2$ the same conclusion holds for every $1\le \alpha
<2^{2/p}$.
\end{corollary}

\begin{proof} As shown in \cite{blmn2} for every $1\le p<\infty$
and $1\le \alpha<2^{\min\{1,2/p\}}$ there is a constant
$c=c(p,\alpha)$ such that for every $k$, $R_p(\alpha,k)\le c\log
k$. It follows from Lemma~\ref{lem:lipcomp} that for every
$\ell\in \N$,
$$
\QS_p^{\Lip}(\alpha,k^\ell)\le k+(c\log
k)\QS_p^{\Lip}(\alpha,k^{\ell-1}).
$$
Since $\QS_p^{\Lip}(\alpha,k)\le k$, by induction we deduce that:
$$
\QS_p^{\Lip}(\alpha,k^\ell)\le \sum_{j=0}^{\ell-1} k(c\log k)^j\le
k(c\log k)^\ell.
$$
Choosing $k$ of the order of $e^{\sqrt{(\log n)(\log\log n)}}$ and
$\ell$ of the order of $\sqrt{\frac{\log n}{\log\log n}}$ gives
the required result.
\end{proof}

We now prove a nearly matching lower bound for $\QS_p(\alpha,n)$.
In order to do so we first observe that Lemma~\ref{lem:aspect}
holds also in the context of Lipschitz quotients.

\begin{lemma}\label{lem:lipaspect} Let $M$ be an $n$-point metric
space and $1<\alpha\le 2$. Then there is a subspace of $M$ which
has an $\alpha$ Lipschitz quotient in an equilateral metric space
and:
$$
|\U|\ge \left\lfloor \frac{n^{(\log \alpha)/[2\log
\Phi(M)]}}{8\log n}\right\rfloor.
$$
\end{lemma}

\begin{proof} The proof is exactly the same as the proof of
Lemma~\ref{lem:aspect}. Using the notation of this proof, the only
difference is that we observe that Lemma~\ref{lem:coloring}
actually ensures that for every $i\neq j$,
$d_M(U_i,U_j),\HH_M(U_i,U_j)\in [\alpha^{\ell-1},\alpha^\ell)$, so
that the quotient obtained is actually a Lipschitz quotient due to
Lemma~\ref{lem:lipq}.

\end{proof}

We also require the following fact, which is essentially proved in
\cite{blmn3} (see Proposition 16 there). Since the result of
\cite{blmn3} was stated for parameters other than what we need
below, we will sketch the proof for the sake of completeness.

\begin{lemma}\label{lem:dichotomy} Fix $0<\e<1$ and let $M$ be an
$n$-point metric space. Then there is a subset $N\subseteq M$
which is either $(1+\e)$-equivalent to an ultrametric, or
$3$-equivalent to an equilateral space, and:
$$
|N|\ge \exp\left(c\sqrt{\frac{\log n}{\log(2/\e)}}\right).
$$
\end{lemma}

\begin{proof} Set $k=2(1/\e+1)$. By Theorem~\ref{thm:kHST} there
is a universal constant $c>0$ such that $M$ contains a subset
$M'\subseteq M$ which is $3$-equivalent to a $(3k)$-HST, $X$, via
a non-contractive bijection $f:M'\to X$, and $|M'|\ge
n^{c/\log(2/\e)}$. Let $T$ be the tree defining $X$. Set
$h=\exp\left(\sqrt{c\frac{\log n}{\log(2/\e)}}\right)$. If $T$ has
a vertex $u$ with out-degree exceeding $h$ then by choosing one
leaf from each subtree emerging from $u$ we obtain a $h$-point
subset of $M'$ which is $3$-equivalent to an ultrametric. We may
therefore assume that all the vertices in $T$ have out-degree at
most $h$. In this case by Lemma 14 in \cite{blmn3} $T$ contains a
binary subtree $S$ with at least $|M'|^{1/\log_2 h}\ge
\exp\left(\sqrt{c\frac{\log n}{\log(2/\e)}}\right)$ leaves. Now,
denote by $\Delta(\cdot)$ the original labels on $S$ (inherited
from $T$). We define new labels $\Delta'(\cdot)$ on $S$ as
follows. For each vertex $u\in S$, denote by $T_1$ and $T_2$ the
subtrees rooted at $u$'s children. We define
$\Delta'(u)=\max\{d_M(x,y);\  x\in f^{-1}(T_1), \, y\in
f^{-1}(T_2)\}$.  As shown in the proof of Case 2 in Proposition 16
of~\cite{blmn3}, this relabelling results in a binary $k$-HST
which is $k/(k-2)=1+\e$ equivalent to $N=f^{-1}(S)$.
\end{proof}

\begin{lemma}\label{lem:lipqlowersmall} For every $1\le \alpha<2$ there is a constant $c=c(\alpha)>0$ such that for every
integer $n$ and every $1\le \alpha<2$,
$$
\QS_2^{\Lip}(\alpha,n)\ge e^{c\sqrt{\log n}}.
$$
\end{lemma}

\begin{proof} By Lemma~\ref{lem:dichotomy} for every $\e>0$ there is a
constant $c=c(\e)$ such that every $n$ point metric space contains
a subset of size at least $e^{c\sqrt{\log n}}$ which is either
$(1+\e)$-equivalent to an ultrametric, in which case we are
already done, or $3$-equivalent to an equilateral space. In the
latter case the subspace obtained has aspect ration at most $3$,
so that the required result follows from
Lemma~\ref{lem:lipaspect}.
\end{proof}

\bigskip

\bigskip
\noindent{\bf Acknowledgements.} We are grateful to Vitali Milman
for proposing the non-linear QS problem to us, and for many
insightful suggestions. We thank Yoav Benyamini and Gideon
Schechtman for helpful comments. We would also like to acknowledge
enjoyable conversions with Shahar Mendel.

\bibliographystyle{plain}

\end{document}